\documentclass[11pt]{amsart}
\usepackage[margin=3cm]{geometry}
\title{Random decompositions of Eulerian statistics}
\author{Alperen \"{O}zdemir}
\address{Department of Mathematics, Georgia Institute of Technology} 
\email{aozdemir6@gatech.edu} 

\usepackage{amssymb, amsmath, float, graphicx, qtree, mathtools, pgfplots, tikz, xcolor}
\usepackage[colorlinks=true, urlcolor=black,linkcolor=blue, citecolor=blue]{hyperref}

\newtheorem{theorem}{Theorem}[section]
\newtheorem{lemma}{Lemma}[section]

\newtheorem{remark}{Remark}[section]

\numberwithin{equation}{section}

\DeclarePairedDelimiter{\floor}{\lfloor}{\rfloor}

\newcommand{\bbox}{\hfill $\Box$}
\newcommand{\pf}{\noindent {\it Proof:} }


\usepackage{tikz}
\usetikzlibrary{arrows,shapes}
\usetikzlibrary{trees}
\usetikzlibrary{matrix,arrows}
\usetikzlibrary{positioning}
\usetikzlibrary{calc,through}
\usetikzlibrary{decorations.pathreplacing}
\usepackage{pgffor}

\usetikzlibrary{decorations.pathmorphing}
\usetikzlibrary{decorations.markings}
\tikzset{
    vector/.style={decorate, decoration={snake}, draw},
	provector/.style={decorate, decoration={snake,amplitude=2.5pt}, draw},
	antivector/.style={decorate, decoration={snake,amplitude=-2.5pt}, draw},
    fermion/.style={draw=black, postaction={decorate},
        decoration={markings,mark=at position .55 with {\arrow[draw=black]{>}}}},
    fermionbar/.style={draw=black, postaction={decorate},
        decoration={markings,mark=at position .55 with {\arrow[draw=black]{<}}}},
    fermionnoarrow/.style={draw=black},
    gluon/.style={decorate, draw=black,
        decoration={coil,amplitude=4pt, segment length=5pt}},
    scalar/.style={dashed,draw=black, postaction={decorate},
        decoration={markings,mark=at position .55 with {\arrow[draw=black]{>}}}},
    scalarbar/.style={dashed,draw=black, postaction={decorate},
        decoration={markings,mark=at position .55 with {\arrow[draw=black]{<}}}},
    scalarnoarrow/.style={dashed,draw=black},
    electron/.style={draw=black, postaction={decorate},
        decoration={markings,mark=at position .55 with {\arrow[draw=black]{>}}}},
	bigvector/.style={decorate, decoration={snake,amplitude=4pt}, draw},
}

\pgfplotsset{compat=1.16}

\begin{document}

\maketitle

\begin{abstract}
This paper develops methods to study the distribution of Eulerian statistics defined by second-order recurrence relations. We define a random process to decompose the statistics over compositions of integers. It is shown that the numbers of descents in random involutions and in random derangements are asymptotically normal with rates of convergence $\mathcal{O} (n^{-1/2})$ and $\mathcal{O}(n^{-1/3})$ respectively.  \end{abstract}

\section{Introduction} \label{intro}
In an earlier work \cite{Oz}, we studied the distribution of the number of descents in random permutations by decomposing the statistic into martingale differences. The technique relied on a submartingale construction by using the first order recurrence relations that the Eulerian numbers satisfy. Then a well-known martingale limit theorem applies and further techniques give rates of convergence in the central limit theorem. The rest of the introduction gives a background to the work and outlines the arguments and results in the paper.

Let $D_n$ be the random variable counting the number of descents in a random permutation of length $n.$ One of the early results on the asymptotic normality of $D_n$ is by Harper's method \cite{Ha}, which applies to statistics whose distribution agrees with the coefficients of a real-rooted polynomial. In particular, $D_n$ can be written as a sum of independent Bernoulli random variables where the probabilities are obtained from the roots of the Eulerian polynomial. Then Lindeberg's condition implies the asymptotic normality, whereas the martingale techniques give the following symmetric decomposition of the same statistic. It was shown in \cite{Oz} that 
\begin{equation*}
D_n=\mathbf{E}(D_n)+\frac{1}{n}\sum_{i=1}^{n} X_i
\end{equation*}
where the distribution of $X_{i+1}$ conditioned on $D_1,\ldots, D_i$ is 
\begin{equation*} 
\begin{aligned}
   &   
   \begin{cases}   
       D_{i}-i,&  \text{with prob. } \frac{D_{i}+1}{i+1},\\[5pt]
      D_{i}+1,&  \text{with prob. } \frac{i-D_{i}}{i+1}.
    \end{cases}
\end{aligned}
\end{equation*}
The martingale techniques also allow us to infer the rates of convergence in the limit theorem matching the rates in classical Berry-Esseen theorems at a cost of calculating the fourth moment compared to the third moment condition of the latter. Nevertheless, the main advantage of the martingale method is that it can be applied to other Eulerian statistics, including those of polynomials that are not real-rooted, with no further difficulty as long as they are defined by a triangular array with a first-order recurrence relation. So we were able to extend our result in various directions, such as to the number of peaks in random permutations, descents in certain conjugacy classes of the symmetric group, descents in Coxeter groups and to vector descent statistics.

This paper addresses Eulerian statistics associated with triangular arrays satisfying second-order recurrence relations. Our first example is the number of descents in random involutions of the symmetric group, which is described in Section \ref{sect2}. Then, in Section \ref{sect3}, we show that it can be written as a martingale difference sequence conditioned on an integer composition. Let $I_n$ be the number of descents in random involutions on $n$ elements. We have
\begin{equation}\label{forminv}
\mathbf{E} \left[I_n\, \big| \,  {C_{\leq 2}(n)=\mathbf{a}}\right]= \mathbf{E}(I_n) + \frac{1}{n}\sum_{i} X_{\mathbf{a_i}}
\end{equation}
where $C_{\leq 2}(n)=(a_1,a_2,\ldots)$ is a random (non-uniform) composition of $a_1+a_2+\ldots=n$ with $a_i=1 \textnormal{ or }2$, $\mathbf{a_i}=(\sum_{j=1}^i a_j,a_i)$ and $\{X_{\mathbf{a_i}}\}_i$ is a martingale difference sequence. The second example is the number of descents in random derangements, which is non-palindromic (asymmetrically distributed about its middle term) unlike $I_n,$ but it still admits a similar decomposition.

In Section \ref{sect3}, in order to derive the decompositions of the statistics and investigate their asymptotic behavior, we define a random process that allows us to decompose them into serially uncorrelated random variables over random compositions of parts of size at most of order of the recurrence relation. Then we study an illustrative example, the number of descents in Fibonacci permutations, and prove various identities that can be of independent interest.
 
In Section \ref{sect4}, we state the limit theorem for martingale difference sequences and prove our main result, which can be combined in the following theorem.
\begin{theorem} \label{main}
Let $I_n$ be the number of descents in random involutions and $R_n$ be the number of descents in random derangements. Then
\begin{align*}
\sup_{x \in \mathbb{R}} \left| \mathbf{P} \left(\frac{I_n - \mathbf{E}[I_n]}{\sqrt{\textnormal{Var}(I_n)}} \leq x\right) - \Phi(x) \right| \hspace*{-1mm} \leq \frac{C}{\sqrt{n}} \quad \textnormal{and} \quad \sup_{x \in \mathbb{R}} \left| \mathbf{P} \left(\frac{R_n - \mathbf{E}[R_n]}{\sqrt{\textnormal{Var}(R_n)}} \leq x\right) - \Phi(x) \right|\hspace*{-1mm} \leq \frac{C'}{\sqrt[3]{n}}
\end{align*}
where $\Phi$ is the standard normal distribution and $C, C'$ are constants  independent of $n$.
\end{theorem}

In Section \ref{sect5}, we provide the extension of the method to recurrences of higher order, which yields martingale difference sequences indexed by compositions with summands of size at most the order of the recurrence relation. Then we give examples of both Eulerian and non-Eulerian statistics defined by second-order recurrence relations that still elude our asymptotic analysis. We conclude with a comment on a possible connection of our work to autoregressive processes.

\section{Eulerian statistics and recurrence relations} \label{sect2}
We first introduce the statistic of interest of the paper, then the  subsets of permutations on which we study it. Let $S_n$ be the symmetric group. A permutation $\pi \in S_n$ is said to have a \textit{descent} at position $i$ if $\pi(i) > \pi(i+1),$ and an \textit{excedance} at position $i$ if $\pi(i)>i.$  A statistic that is equidistributed with the number of descents or excedances is called an \textit{Eulerian statistic}. In fact, the number of descents is equidistributed with the number of excedances over uniformly random permutations.  

The number of permutations with a given number of descents are counted by \textit{Eulerian numbers}. Let $A_{n,k}$ be the number of permutations of $n$ elements with $k$ descents. We define
\begin{equation*}
A_n(t)=\sum_{\pi \in S_n}t^{\textnormal{des}(\pi)}=\sum_{k\geq 0} A_{n,k} t^k,
\end{equation*}
where $\textnormal{des}(\pi)$ is the number of descents in $\pi \in S_n.$ $A_n(t)$ is called the \textit{Eulerian polynomial} and has the rational generating function
\begin{equation*}
\frac{A_n(t)}{(1-t)^{n+1}}= \sum_{k \geq 0} (k+1)^n t^k.
\end{equation*}
Elementary manipulations in the sum on the right-hand side will give the following recursive relation on the Eulerian numbers.
\begin{equation*}
A_{n+1,k} = (k+1) A_{n,k} + (n-k+1) A_{n, k-1}.
\end{equation*}
The transition probabilities for the martingale associated with $D_n$ are obtained from these recursive relations, which is elaborated in Section 5 of \cite{Oz} with examples in Section 6 of the same article. We will use those techniques in the following two examples. 

\subsection{Involutions}

 A permutation $\pi$ is called an \textit{involution} if its inverse is itself; in other words, $\pi^2=1.$ Involutions consist of combinations of fixed points and transpositions. They do not form a conjugacy class, but can be thought of as the union of conjugacy classes of permutations that only consist of fixed points and $2-$cycles. Before turning to the number of descents in random involutions, we first consider the involutions themselves.

Let $i_n$ denote the number of involutions of length $n.$ We observe that the number of involutions of length $n$ where $n$ is a fixed point is $i_{n-1}$, and the number of involutions where $(i  \, n)$ is a cycle is $i_{n-2}$ for $1 \leq i \leq n-1.$ Therefore we have the following second-order recursive relation
\begin{equation} \label{tele}
 i_n= i_{n-1}+ (n-1) i_{n-2},
 \end{equation} 
which is employed by Rothe to list $i_1$ to $i_{10}$ \cite{H00}. A comphrehensive account on involutions can be found in \cite{K3}, which includes the asymptotic formula 
\begin{equation*}
i_n = \frac{1}{\sqrt[4]{4e}} \left(\frac{n}{e}\right)^{n/2}e^{\sqrt{n}}\left(1 + \mathcal{O}\left(n^{-1/2}\right)\right).
\end{equation*}
It was shown earlier in \cite{CHM51} that the growth rate $\frac{i_{n+1}}{i_{n}}$ is asymptotically of order $\sqrt{n}$, more precisely 
\begin{equation}\label{invratio}
 \sqrt{n+1}\leq \frac{i_{n+1}}{i_n} \leq \sqrt{n+1} +1.
 \end{equation} 
 
An interesting fact is that the number of involutions in the symmetric group $S_n$ is equal to the number of \textit{standard Young tableaux} of size $n.$ A standard Young tableau of size $n$ is a diagram with $n$ boxes of left-justified rows of non-increasing lengths, such that each box is filled with a number in increasing order from both left to right and top to bottom. Each irreducible representation of $S_n$ is represented by a diagram and the number of ways to fill it in the described way gives the dimension of the representation, 
\begin{equation}\label{invdim}
i_n= \sum_{\lambda \textnormal{ irr.}} \textnormal{dim}(\lambda).
\end{equation}
See, for instance, \cite{S13}. In fact, the dimensions of irreducible representations of any group add up to the number of involutions if Frobenius-Schur indicators of its characters are equal to $1$, which is the case for Coxeter groups of type A (symmetric groups), B and D \cite{GP00}. We do not examine this topic further here except two identities \eqref{stan1} and \eqref{stan2} in the following section, but note that our results are easily extendible to Coxeter groups of type $B.$ 
\begin{figure}[H]
  \begin{tabular}{rccccccccccc}
$n=1$:\quad & &    &    &    &    &  1\\\noalign{\smallskip\smallskip}
$n=2$:\quad &  &   &    &    & 1  &    &  1\\\noalign{\smallskip\smallskip}
$n=3$:\quad &  &  &    & 1  &    &  2 &    &  1\\\noalign{\smallskip\smallskip}
$n=4$:\quad &  &  &  1 &    & 4  &    &  4 &    &  1\\\noalign{\smallskip\smallskip}
$n=5$:\quad & & 1 &    &  6 &    &  12 &    &  6 &    &  1\\\noalign{\smallskip\smallskip}
$n=6$:\quad &  1 &    &  9 &    &  28 &    &  28 &    &  9 & & 1\\\noalign{\smallskip\smallskip}
\end{tabular}
\vspace*{-0.4 cm}
\caption{The number of involutions of length $n$ with $k$ descents, $I_{n,k}$}
\end{figure}

The generating function for the number of involutions with a given number of descents is obtained in \cite{DF} by identities involving Schur functions and Robinson-Schensted correspondence between involutions and standard Young tableaux. It is also derived in \cite{GR93} from quasisymmetric functions of descent compositions as a result of a more general formula including major indices. It satisfies the folowing rational form.
\begin{equation*}
\sum_{n \geq 0}I_n(t)\frac{u^n}{(1-t)^{n+1}} = \sum_{k\geq 0} \frac{t^k}{(1-u)^{k+1}(1-u^2)^{\binom{k+1}{2}}},
\end{equation*}
where $$I_{n}(t)=\sum_{\{\pi: \pi^2=1 \}} t^{\textnormal{des}(\pi)} = \sum_{k=1}^{n-1}I_{n,k}t^k.$$
The coefficient of $u^n$ is given by
\begin{equation} \label{invoconvo}
\frac{I_n(t)}{(1-t)^{n+1}}= \sum_{r=0}^{\infty} \sum_{s=0}^{\floor{n/2}}\binom{\binom{r+1}{2}+s-1}{s}\binom{n-r-2s}{n-2s}t^r.
\end{equation}

It is further studied in \cite{GZ} to show the unimodalty of $I_n(t).$ It fails to be a log-concave sequence, see \cite{Ba}. The gamma-positivity of the polynomial is proved in \cite{W19}. Applying Zeilberger's algorithm to the coefficient of $t^r$ in \eqref{invoconvo}, the following formula is obtained in \cite{GZ}. For $n\geq 1,$

\begin{align}\label{invosecond}
\begin{split}
nI_{n}(t)
&=(t-t^2)I_{n-1}'(t)+[1+(n-1)t]I_{n-1}(t)
+t^2(1-t)^2I_{n-2}''(t) \\[5pt]
&\quad+t(1-t)[3+(2n-5)t]I_{n-2}'(t)
+(n-1)[1+t+(n-2)t^2]I_{n-2}(t).
\end{split}
\end{align}
Observe that $I_{n}(1)$ is equal to the number of involutions, $i_n,$ and \eqref{tele} can also be verified by the formula above. Comparing the coefficients in \eqref{invosecond}, the second-order recurrence relation below is shown in \cite{GZ}.
\begin{align}\label{coeff}
\begin{split}
I_{n+2,k}&= \frac{k+1}{n+2} I_{n+1,k} + \frac{n-k+2}{n+2} I_{n+1,k-1} \\
&+ \frac{(k+1)^2 +n}{n+2} I_{n,k} + \frac{2k(n-k+1)-n+1}{n+2} I_{n,k-1}+\frac{(n-k+2)^2 +n}{n+2} I_{n,k-2}.
\end{split}
\end{align}
By the symmetry of the coefficients, $I_{n+2,k}=I_{n+2,n+1-k}.$ Therefore $I_n(t)$ is a palindromic polynomial, which implies  $\mathbf{E}(I_n)=\frac{n-1}{2}.$

Next, we invoke the method in \cite{Oz} to obtain martingales from the recursive relations. Let $I_n$ be the random variable counting the number of descents in random involutions. Observe that $I_1=0.$ Fixing $k,$ we first write down the recursive expansions of $I_{n+2,k+1}$ and $I_{n+2,k+2}$ in addition to \eqref{coeff} as they include the terms $I_{n,k}$ or $I_{n+1,k},$ whose coefficients will be the transition probabilities. 
\begin{align*}
I_{n+2,k+1}&= \ldots + \frac{n-k+1}{n+2} I_{n+1,k}+\ldots+\frac{2(k+1)(n-k)-n+1}{n+2} I_{n,k} \\
I_{n+2,k+2}&=\ldots + \frac{(n-k)^2 +n}{n+2} I_{n,k}. 
\end{align*}
We replace $k$ by $I_n$ and $I_{n+1}$ in the coefficients of $I_{n,k}$ and $I_{n+1,k}$ respectively to obtain the following process.
\begin{equation} \label{invjump}
I_{n+2}=\begin{aligned}
   &   
   \begin{cases}
     I_{n},&  \text{with prob. } \frac{(n+1)i_n}{i_{n+2}}\frac{(I_n+1)^2+n}{(n+1)(n+2)},\\[5pt]
      I_{n}+1,&  \text{with prob. } \frac{(n+1)i_n}{i_{n+2}}\frac{2(I_n+1)(n-I_n)-n+1}{(n+1)(n+2)},\\[5pt]
       I_{n}+2,&  \text{with prob. } \frac{(n+1)i_n}{i_{n+2}}\frac{(n-I_n)^2+n}{(n+1)(n+2)},\\[5pt]
       I_{n+1},&  \text{with prob. } \frac{i_{n+1}}{i_{n+2}}\frac{I_{n+1}+1}{n+2},\\[5pt]
      I_{n+1}+1,&  \text{with prob. } \frac{i_{n+1}}{i_{n+2}}\frac{n-I_{n+1}+1}{n+2}.
    \end{cases}
\end{aligned}
\end{equation}

\subsection{Derangements}\label{sect:der}

A \textit{derangement} $\pi$ is a permutation with no fixed points; in other words $\pi(i)\neq i$ for all $i.$ A standard method to count the derangements is the inclusion-exclusion principle. See \cite{C11} for two other methods and a short algebraic survey on them. The number of derangements is defined recursively by
\begin{equation}\label{dersecrec}
d_n=n d_{n-1} +(-1)^n=(n-1)(d_{n-1}+d_{n-2}).
\end{equation}
It has the closed-form expression
\begin{equation}\label{derclose}
 d_n=n! \sum_{i=0}^n \frac{(-1)^i}{i!},
\end{equation}
which is asymptotically $\frac{n!}{e}.$ In fact, $d_n$ is the closest integer to $\frac{n!}{e}$ since the error term in the Taylor expansion of $e^t$ evaluated at $t=-1$ is factorially small. The generating function for the number of descents in derangements is 
\begin{equation*}
 \sum_{n\geq 0} D_n(t) \frac{z^n}{(1-t)^{n+1}} = \sum_{k \geq 1} \frac{(1-z)^kt^k}{(1-kz)}
  \end{equation*}
  where
  \begin{equation*}
  D_{n}(t)=\sum_{\{\pi:\, \forall i \,\pi(i)\neq i  \}} t^{\textnormal{des}(\pi)} = \sum_{k=1}^{n-1}D_{n,k}t^k.
  \end{equation*}
The generating function is derived in the same paper of Gessel and Reutenauer \cite{GR93} mentioned in the preceding section. Evaluating the coefficient of $z^n$ on the right-hand side, we have 
\begin{equation}\label{deragen}
  \frac{D_n(t)}{(1-t)^{n+1}} = \sum_{k\geq 0} \left(\sum_{i=0}^n (-1)^i  \binom{k}{i} k^{n-i} \right) t^k
  \end{equation}  
 The above expression is studied in \cite{FZ18}, and they obtained
 \begin{equation} \label{deralt}
 D_n(t)=(-1)^n t^{n-1} + (1+(n-1)t)D_{n-1}(t) +t(1-t)D'_{n-1}(t),
 \end{equation}
for $n\geq 2.$ Although $D_n(t)$ is not palindromic as displayed in Figure \ref{figd}, it is unimodal and the maximum coefficient appears in the middle, which is shown in \cite{FZ18}. The same properties were obtained earlier in \cite{Z96} for the number of excedances in derangements.
 \begin{figure}[H]
  \begin{tabular}{rccccccccccc}
$n=2$:\quad \quad & &    &    &    &    &  1\\\noalign{\smallskip\smallskip}
$n=3$:\quad \quad &  &   &    &    & 2  &    &  0\\\noalign{\smallskip\smallskip}
$n=4$:\quad \quad &  &  &    & 4  &    &  4 &    &  1\\\noalign{\smallskip\smallskip}
$n=5$:\quad \quad &  &  &  8 &    & 24  &    &  12 &    &  0\\\noalign{\smallskip\smallskip}
$n=6$:\quad \quad & & 16 &    &  104 &    &  120 &    &  24 &    &  1\\\noalign{\smallskip\smallskip}
$n=7$:\quad \quad &  32 &    &  392 &    &  896 &    &  480 &    &  54 & & 0\\\noalign{\smallskip\smallskip}
\end{tabular}
\vspace*{-0.4 cm}
\caption{The number of derangements of length $n$ with $k\geq 1$ descents, $D_{n,k}$}
\label{figd}
\end{figure}

In order to eliminate the alternating term, we expand $tD_{n-1}(t)$ according to  \eqref{deralt} to have
\begin{align*}
D_{n}(t)=(1+(n-2)t)D_{n-1}(t) +t(1-t)D'_{n-1}(t) + (1+(n-2)t)D_{n-2}(t) +t(1-t)D'_{n-2}(t).
\end{align*}
Equating the coefficients on both sides, we obtain
\begin{equation}
D_{n,k}=(k+1)D_{n-1,k}+(n-k-1)D_{n-1,k-1} + k D_{n-2,k-1} + (n-k)D_{n-2,k-2}.
\end{equation}
As in the case with the involutions, we define the following random process that gives the number of descents in random derangements of permutations of length $n$ at its $n$th stage. We denote it by $R_n$ in order to avoid confusion with the number of descents in random permutations, and note that $R_1=0.$ The random sequence satisfies
\begin{equation} \label{derjump}
R_{n+2}=\begin{aligned}
   &   
   \begin{cases}
     R_{n}+1,&  \text{with prob. } \frac{(n+1)d_n}{d_{n+2}}\frac{R_n+1}{(n+1)},\\[5pt]
     R_{n}+2,&  \text{with prob. } \frac{(n+1)d_n}{d_{n+2}}\frac{n-R_n}{(n+1)},\\[5pt]
       R_{n+1},&  \text{with prob. } \frac{(n+1)d_{n+1}}{d_{n+2}}\frac{R_{n+1}+1}{(n+1)},\\[5pt]
      R_{n+1}+1,&  \text{with prob. } \frac{(n+1)d_{n+1}}{d_{n+2}}\frac{n-R_{n+1}}{(n+1)}
    \end{cases}
\end{aligned}
\end{equation}

Having defined $I_n$ and $R_n$ as random processes, we are able to treat them in a general setting.

\section{Random decompositions}\label{sect3}
The aim of this section is twofold. It is to lay the groundwork for the study of the statistics described in Section \ref{sect2} and to deliver the intuition for the recursive methods of the next section by an example. Regarding the former, here and in Section \ref{sect:bincomb}, we show how to transform the recursive relations defined by second-order recurrences as in \eqref{invjump} and \eqref{derjump} to a forward moving process. This can be compared to the martingale derivation for descent statistics from the first-order recurrences in Chapter 5 of \cite{Oz}. In the second part, Section \ref{sect:fib} and \ref{sect:app}, we give a simple example for the process, and develop an approximation method by demonstrating its use on the example, which is then to be used in Chapter 4. We also derive some combinatorial identities of independent interest. 

Consider a stochastic process $\{Z_n\}_{n \geq 1}$ obeying the following rule.  
\begin{equation} \label{nonmarkov}
Z_{n}=\begin{aligned}
   &   
   \begin{cases}
    Z_{n-2}+f(Z_{n-2}),&  \text{with prob. } q_{n},\\
    Z_{n-1}+f(Z_{n-1}), &  \text{with prob. } 1-q_{n},\\
    \end{cases}
\end{aligned}
\end{equation}
where $f$ is some measurable function. This process can be viewed as a general form of \eqref{invjump} and \eqref{derjump}. The martingale methods do not immediately apply to $\{Z_n\}_{n \geq 1},$ since it is non-Markovian,
\begin{equation*}
\mathbf{E}[Z_{n}|\mathcal{F}_{n-1}]= (1-q_n) Z_{n-1} + q_n Z_{n-2}
\end{equation*}
where $\mathcal{F}_i$ is the $\sigma-$field generated by $Z_1,\dots,Z_i$ and assuming that the increments $f(Z_{n-2})$ and $f(Z_{n-1})$ have zero mean. The process can be updated from $Z_{n-2}$ besides $Z_{n-1}$; we call an update a $two-jump$ in the former case and a $one-jump$ in the latter. 
\subsection{Binary words to compositions} \label{sect:bincomb}

We now describe how to split the random variable $Z_n$ once we keep track of the jumps of $Z_n$ as whether it is updated from $Z_{n-2}$ by a two-jump with probability $q_n$ or from $Z_{n-1}$ by a one-jump with probability $1-q_{n}$. We assume that the first update is always a one-jump. What we end up with is a binary word $w$ of $n$ letters which starts with $1$. This does not directly give us the decomposition of $Z_n$. But if we discard all entries left to 2s starting from the rightmost letter of the word by disregarding the already discarded 2s, we can label the indices of the decomposition by the remaining ones. Let $\psi$ denote the mapping of the discard operation, for example
\vspace{-5mm}
\begin{center}
\begin{align} \label{onetwoone}
\begin{split}
w: \quad \mathbf{1} \quad & 2 \quad 1 \quad 2 \quad 2 \quad 1 \quad 2 \quad 1 \quad 1 \quad 2 \quad 2 \quad 2 \\ 
\psi(w): \quad \, \, \, \phantom{} \quad & 2 \quad 1 \quad \phantom{0} \quad 2 \quad \phantom{0} \quad 2 \quad 1 \quad \phantom{0} \quad 2 \quad \phantom{0} \quad 2
\end{split}
\end{align}
\end{center}
The latter string is indeed a composition of $n$ of summands of size at most $2.$ We denote the set of compositions of $n$ consisting of only $1$ and $2$ by $\mathcal{C}_{\leq 2}(n)$.

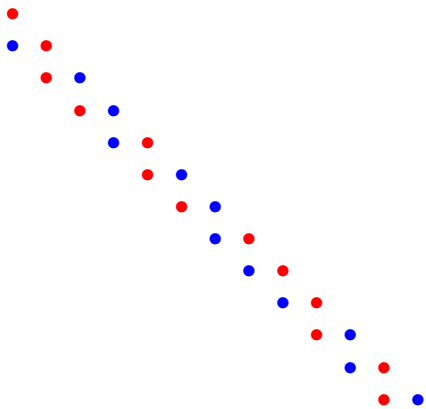
\begin{figure}[h!]
\begin{minipage}{0.5\linewidth}
\begin{tikzpicture}
\matrix [matrix of math nodes]
{
 |[red]| \bullet &   \\
|[blue]| \bullet  & |[red]| \bullet &   \\
 &  |[red]| \bullet & |[blue]| \bullet & & \\
 & & |[red]| \bullet & |[blue]| \bullet & & \\
& & & |[blue]| \bullet & |[red]| \bullet & \\
& & & &|[red]| \bullet & |[blue]| \bullet \\
& & & & &|[red]| \bullet & |[blue]| \bullet \\
& & & & & &|[blue]| \bullet & |[red]| \bullet \\
& & & & & & &|[blue]| \bullet & |[red]| \bullet \\
& & & & & & & &|[blue]| \bullet & |[red]| \bullet \\
& & & & & & & & & |[red]| \bullet & |[blue]| \bullet \\
& & & & & & & & & &|[blue]| \bullet & |[red]| \bullet \\
& & & & & & & & & & &|[red]| \bullet & |[blue]| \bullet \\
};
\end{tikzpicture}
\end{minipage}
\begin{minipage}{0.45\linewidth}
Starting from the bottom, we move along the diagonal path on the right in the north-west direction to obtain the decomposition. If the color alters, it means that process is updated by a two-jump, $X_{*,2},$ otherwise it is updated by $X_{*,1}.$ For instance, if we take the unpaired particle $Z_0=0$ and update the positions of the others according to the rule \eqref{rule}, then we end up with $Z_{12}$ which decomposes as $Z_{12}=\mathbin{\color{blue}{X_{12,2}}}  +  \mathbin{\color{blue}{X_{10,2}}} + \mathbin{\color{red}{X_{8,1}}}+  \mathbin{\color{red}{X_{7,2}}}  + \mathbin{\color{blue}{X_{5,2}}} +  \mathbin{\color{blue}{X_{3,1}}} + \mathbin{\color{blue}{X_{2,2}}}$ where the color indicates the particle on which the update is based.
  \vfill
\end{minipage}
  \caption{An example of an update rule to obtain the decomposition of $\{Z_n\}_{n \geq 1}.$ The rows of particles in the figure represent $\zeta_0, \zeta_1,\ldots,\zeta_{12}.$}
  \label{fig1}
\end{figure}

Next, we formalize the idea by defining a process that generates integer compositions, which will allow us to decompose $Z_n$ over its outcomes and express it in the martingale setting. We first generate a step-ahead copy of the process and then couple two processes in the following way to retain the information from the previous stage and to define the two types of jumps at every stage. Consider a vector $\zeta_n=\left(Z_n,Z_{n-1}\right)$ of random variables with the update rule
\begin{equation}\label{rule}
\zeta_{n+1}=\begin{cases}
\left(Z_{n-1}+X_{n+1,2}, Z_n\right) \textnormal{ with prob. } q_{n+1}, \\[5pt]
\left(Z_n+X_{n+1,1}, Z_n\right) \textnormal{ with prob. } 1-q_{n+1} 
\end{cases}
\end{equation}
where $X_{n+1,1}$ and $X_{n+1,2}$ may depend on $\zeta_n.$ We take $\zeta_0=Z_0.$ One-jump means that the particles move one step forward together, while two-jump means that the particle in the back moves two steps forward and the particle in the front stays at its position. See Figure \ref{fig1} and the accompanying explanation below. 

We denote the composition obtained from $\zeta_n$ by the random variable $C_{\leq 2}(n)$, which is recorded by the secondary subindices of differences. For the example above, 
\begin{equation*}
\mathbf{P}(C_{\leq 2}(12)=(2,1,2,2,1,2,2))=q_2 (1-q_3) q_5 q_7 (1-q_8) q_{10} q_{12}.
\end{equation*}

We color the summands in Figure \ref{fig1} to indicate the dependence of the jumps on the preceding stage. However, if the jumps form a martingale difference sequence, i.e., $\mathbf{E}[X_{i,1}|\mathcal{F}_{i-1}]=\mathbf{E}[X_{i,2}|\mathcal{F}_{i-1}]=0$ for all $i,$ then the random variable from which the difference is obtained has no correlation with it, thus no color is needed. In Section \ref{sect4}, either the jumps are already in that form (in the case of involutions), or will be shown to be close enough to that form (in the case of derangements) for the asymptotic analysis. 
 
\subsection{A deterministic example and some identities} \label{sect:fib}

In order to explain the idea in the previous section with a simple example, for which $f$ in \eqref{nonmarkov} is deterministic, we consider the \textit{Fibonacci permutations}. A permutation $\pi$ is called a Fibonacci permutation if $|\pi(i)-i| \leq 1.$ The restriction implies that $\pi$ consists only of fixed points and pairwise adjacent transpositions, which makes it an involution. For instance,
\begin{equation*}
\pi= 1\,(2 \,3) \,  4 \, 5 \, (6 \,7)\, 8 \, (9 \, 10)
\end{equation*}
is a Fibonacci permutation. If we take $2$s as transpositions and $1$s as fixed points in the second line of \eqref{onetwoone}, we have a bijection between Fibonacci permutations of length $n$ and $\mathcal{C}_{\leq 2}(n).$
They are counted by the $(n-2)$nd Fibonacci number. Fibonacci numbers satisfy the simplest second-order recurrence relation,
\begin{equation*}
f_{n}=f_{n-1}+f_{n-2}
\end{equation*}
with the initial condition that $f_0=f_1=1.$ This recursion can be verified for the number of Fibonacci permutations by the process of either adding $n$ as a fixed point to a permutation of size $n-1$ or adding the transposition $\left(n \, \, n-1\right)$ to a permutation of size $n-2$.

Let $F_n$ be the number of descents in random Fibonacci permutations. Observe that the number of descents agrees with the number of transpositions in these permutations. Moreover, the short proof of the recursion above defines a descent counting process; the former case has no contribution to the number of descents while the latter contributes by $1$. This gives
\begin{equation}\label{fibos}
F_{n+2}=\begin{aligned}
   &   
   \begin{cases}
     F_{n}+1,&  \text{with prob. } \frac{f_n}{f_{n+2}},\\[5pt]
      F_{n+1},&  \text{with prob. } \frac{f_{n+1}}{f_{n+2}}.
    \end{cases}
\end{aligned}
\end{equation}
Each jump takes value either $0$ or $1$ with probabilities independent of $F_n.$

We can directly compute the mass function of $F_{n+2}$ by using the map in \eqref{onetwoone}. By \eqref{fibos}, $F_{n+2}$ is equal to $k$ if and only if there are $k$ gaps, associated with $k$ two-jumps, in the second line of \eqref{onetwoone}. The following sum for $\mathbf{P}(F_{n+2}=k)$ is over all possible locations of gaps, denoted by $(j_1,\ldots,j_k).$ We also need to consider the all possibilities for discarded jumps in the first line for $w$ in \eqref{onetwoone}. There are two possibilities for each gap, either a one-jump or a two-jump. In the first line of the sum below, these probabilities are respectively $\frac{f_{j-2}}{f_{j-1}}$ and $\frac{f_{j-3}}{f_{j-1}}$ which immediately precede the probability of the two-jump, $\frac{f_{j-2}}{f_j}$. 
\begin{align*}
\mathbf{P}(F_{n+2}=k)=&\sum_{1<j_1 <\cdots < j <\cdots < j_k}   \frac{\cdots}{\cdots f_{j-2}} \left(\frac{f_{j-2}}{f_{j-1}} \frac{f_{j-2}}{f_{j}} + \frac{f_{j-3}}{f_{j-1}} \frac{f_{j-2}}{f_{j}}\right) \frac{\cdots}{f_{j+1}\cdots}  \\
=& \sum_{1<j_1 <\cdots < j <\cdots < j_k} \frac{f_1 \cdots (f_{j-3}+f_{j-2})f_{j-2}\cdots f_{n+1}}{f_2 \cdots f_{j-2}f_{j-1}f_{j}f_{j+1}\cdots f_{n+2}} \\
=&\frac{1}{f_{n+2}}\sum_{1<j_1 <\cdots < j_k} 1\\
=&\frac{\binom{n-k}{k}}{f_{n+2}}.   
\end{align*}

See \cite{DGH01} for a shorter argument for the above result and the use of the following generating function 
\begin{align*}
f(t)=\sum_{k=0}^{\floor{n/2}} \binom{n-k}{k}t^k = \frac{1}{\sqrt{1+4t}}\left(\left( \frac{1+ \sqrt{1+4t}}{2}\right)^{n+1} - \left( \frac{1- \sqrt{1+4t}}{2}\right)^{n+1} \right)
\end{align*}
to calculate
\begin{equation} \label{fibomom}
\mathbf{E}(F_n) = \frac{5-\sqrt{5}}{10} n + \frac{1-\sqrt{5}}{10} + \mathcal{O}(e^{-n}) \, \textnormal{  and  } \, \textnormal{Var} (F_n)=\frac{n}{5 \sqrt{5}}+ \mathcal{O}(1).
\end{equation}
The asymptotic normality can be obtained by the methods listed in \cite{DGH01}. A particularly interesting one is an application of Harper's method to the number of edges in random matchings of graphs, where $F_n$ is associated with $(n-1)$-path \cite{G81}. Furthermore, observing that the number of edges in the matchings of the complete graph $K_n$ on $n$ vertices is equidistributed with the number of transpositions in random involutions, the author of \cite{G81} shows the asymptotic normality of the latter as a corollary.

Next, we apply the same idea for the computation of $\mathbf{P}(F_{n+2}=k)$ to obtain identities from the recurrence relations that the number of involutions and derangements obey. Considering the jump probabilities $\frac{(n+1)i_n}{i_{n+2}}$ and $\frac{i_{n+1}}{i_{n+2}}$ in the case of involutions, we sum over the number of 2-jumps (or, equivalently, the number of gaps in \eqref{onetwoone}), which we denote by $k,$ to have

\begin{align*}
1\hspace*{-0.5mm}=& \hspace*{-1mm} \sum_{k=1}^{\floor{\frac{n+2}{2}}} \sum_{1<j_1 <\cdots < j <\cdots < j_k} \hspace*{-4mm} \frac{\cdots}{\cdots  i_{j-2}} \left( \frac{i_{j-2}}{i_{j-1}} \frac{(j-1)i_{j-2}}{i_{j}} +\frac{(j-2)i_{j-3}}{i_{j-1}} \frac{(j-1)i_{j-2}}{i_{j}}\right) \frac{\cdots}{i_{j+1} \cdots}  \\
=& \hspace*{-1mm} \sum_{k=1}^{\floor{\frac{n+2}{2}}}\sum_{1<j_1 <\cdots < j <\cdots < j_k} \prod_{s=1}^k (j_s-1) \frac{i_1 \cdots [(j-2)i_{j-3}+i_{j-2}]i_{j-2}\cdots i_{n+1}}{i_2 \cdots i_{j-2}i_{j-1}i_{j}i_{j+1}\cdots i_{n+2}} \\
=&\frac{1}{i_{n+2}}\sum_{k=1}^{\floor{\frac{n+2}{2}}} \sum_{1<j_1 <\cdots < j <\cdots < j_k} \prod_{s=1}^k (j_s-1).  
\end{align*}
Thus, by \eqref{invdim}, we obtain
\begin{equation}\label{stan1}
 \sum_{\mathbf{a} \in \mathcal{C}_{\leq 2}(n)}\prod_{\{i : a_i=2\}} (1+a_1+\cdots+a_{i-1})= i_{n+1}=\sum_{\lambda  \vdash n+1}\textnormal{dim}(\lambda)
\end{equation}
where $\lambda$ runs over the partitions of $n+1,$ which in fact label the irreducible representations of the symmetric group. If we take $n!$ instead of $i_n$ above and observe that  $n!=(n-1)!+(n-1)^2 (n-2)!,$ we have
\begin{equation}\label{stan2}
 \sum_{\mathbf{a} \in \mathcal{C}_{\leq 2}(n)}\prod_{\{i : a_i=2\}} (1+a_1+\cdots+a_{i-1})^2=(n+1)!= \sum_{\lambda  \vdash n+1}\textnormal{dim}(\lambda)^2
\end{equation}
by the well-known fact in representation theory that the sum of the squares of the dimensions of irreducible representations of a group is equal to its size.
The formulas \eqref{stan1} and \eqref{stan2} (communicated to Persi Diaconis by Richard Stanley \cite{D18}) in terms of bits of binary strings with no two consecutive ones are already stated in \cite{D18}. 
 
Similarly, using \eqref{dersecrec} in the case of derangements,
\begin{align*}
\begin{split} 
1=& \sum_{k=1}^{\floor{\frac{n+2}{2}}} \sum_{2<j_1 <\cdots < j <\cdots < j_k}  \frac{\cdots}{\cdots d_{j-2}}\left( \frac{(j-2)d_{j-2}}{d_{j-1}} \frac{(j-1)d_{j-2}}{d_{j}}  + \frac{(j-2)d_{j-3}}{d_{j-1}} \frac{(j-1)d_{j-2}}{d_{j}} \right) \frac{\cdots}{d_{j+1}\cdots} 
\end{split}
\\
=& \sum_{k=1}^{\floor{\frac{n+2}{2}}}\sum_{2<j_1 <\cdots < j <\cdots < j_k}  \frac{d_1 \cdots [(j-2)(d_{j-3}+d_{j-2})](j-1)d_{j-2}j\cdots (n+1) d_{n+1}}{d_2 \cdots d_{j-2}d_{j-1}d_{j}d_{j+1}\cdots d_{n+2}} \\
=&\frac{(n+1)!}{d_{n+2}}\sum_{k=1}^{\floor{\frac{n+2}{2}}} \sum_{2<j_1 <\cdots < j <\cdots < j_k} \prod_{s=1}^k \frac{1}{j_s-2}  .
\end{align*}
We then obtain the following identity by \eqref{derclose}.
\begin{equation*}
\frac{1}{n+2} \sum_{\mathbf{a} \in \mathcal{C}_{\leq 2}(n)}\prod_{\{i : a_i=2\}} \frac{1}{a_1 +\cdots+a_i}=\sum_{k=0}^{n+2} \frac{(-1)^k}{k!}, 
\end{equation*}
recalling that $\mathbf{a}=(a_1, a_2,\ldots)$ is a composition of $n$ with parts of size $1$ or $2$.

\subsection{Approximated distributions}\label{sect:app}

We have discussed how the randomness is inherited from binary words to compositions over several examples. The exact expressions for statistics over compositions are usually difficult to derive. In the remainder of this section, we discuss how to approximate those expressions with the help of binary word statistics. We start with the Fibonacci permutations. The ratio $\frac{f_{n+1}}{f_n},$ which defines the probabilities in \eqref{fibos}, is well-known to converge to $\varphi:=\frac{1+\sqrt{5}}{2}$ and satisfy the equation $\varphi^2=1+ \varphi.$ The error term in the ratio is   
\begin{equation}\label{appfibo}
\left|\frac{f_{n+1}}{f_{n}} - \varphi\right| \leq \frac{1}{\varphi^n f_{n+2}},
\end{equation}
see Section 1.35 of \cite{V12}.
The expected value of $F_n$ can be estimated by the expected number of two-jumps minus the expected number of adjacent two-jumps since only two-jumps contribute to $F_n$. Then to compensate for the case with three two-jumps in a row, we add its expected value and so on. This will give us an inclusion-exclusion argument. Considering that the probability of a two-jump is approximately $\frac{1}{\varphi^2}$ by \eqref{appfibo}, we have
\begin{equation*}
\mathbf{E}(F_n) \approx \sum_{i=0}^n \frac{(-1)^i(n-i)}{\varphi^{2i+2}}=\frac{5-\sqrt{5}}{10} n +\mathcal{O}(1).
\end{equation*}

In general, we consider a binary random variable $\mathcal{B}(p,\alpha, \beta)$ that is equal to $\alpha$ with probability $p$ and equal to $\beta$ otherwise. If the mean of $\mathcal{B}(p,\alpha, \beta)$ is zero, we will write it as $\mathcal{B}(p,\alpha).$  Let $T_n=\sum_{k=1}^n \mathcal{B}(q_k, \beta_{k,1}, \beta_{k,2}).$ Now, we associate the summands of $T_n$ with letters of a binary word and define a new statistic. Let $w$ be a binary word of length $n.$ By an abuse of notation, we define the random variable $\psi(T_n)$ induced by the function in \eqref{onetwoone}, 
\begin{equation}
\psi (T_n) (w)= \sum_{i} \beta_{\mathbf{a_i}}
\end{equation}
where $\psi(w)=(a_1, a_2,\dots)$ and $\mathbf{a_i}=\left(\sum_{j=1}^i a_j,a_i\right)$. For instance, if we take $T_n$ to be $\sum_{i=2}^n \mathcal{B}\left(f_{i-2}/f_{i},0,1\right)$, then $\psi(T_n)=F_n.$ 
 Although we will derive more indirect bounds for $\psi(T_n)$ in our two main examples, we provide the following upper bound on the variance of $\psi(T_n)$ that could be useful in the study of different examples.
\begin{lemma}\label{lemma}
 Let $T_n=\sum_{i=1}^n \mathcal{B}(p_i,\beta_i)$ where the summands are independent zero-mean random variables. Then
$\textnormal{Var}(\psi(T_n))\leq \textnormal{Var}(T_n).$
\end{lemma}
\pf By the conditional variance formula,
 \begin{equation}\label{varfor}
 \textnormal{Var}(T_n)= \mathbf{E}[\text{Var}(T_n| \, C_{\leq 2}(n))]+
 \text{Var}(\mathbf{E}[T_n| \, C_{\leq 2}(n)]).
 \end{equation}
Since the summands of $T_n$ have zero mean, 
\begin{equation*}
\mathbf{E}[T_n | \, C_{\leq 2}(n)]=\mathbf{E}[\psi(T_n) |\, C_{\leq 2}(n)].
\end{equation*}
Observing that $\psi(T_n)$ is deterministic conditioned on the composition $C_{\leq 2}(n),$ 
\begin{equation*}
\textnormal{Var} (\mathbf{E}[T_n |\, C_{\leq 2}(n)])= \textnormal{Var}(\psi(T_n))
\end{equation*}
The result follows from \eqref{varfor}.
\bbox

\begin{remark}
If we take $T_n$ to be the sum of Bernoulli random variables with an identical success probability and $1$ assigned to two-jumps, such as in the approximation to $\mathbf{E}(F_n)$ with $p=1/\varphi,$ the conclusion of Lemma \ref{lemma} holds true, which can be shown by conditioning on the number of two-jumps. 
\end{remark}

\section{Limit theorem and the rate of convergence}\label{sect4}

In this section, we study the asymptotic distributions of the processes defined in Section \ref{sect2} by martingale limit theorems. We already mentioned that these processes do not yield martingales; however, we can write them as random sums of martingale differences as outlined in Section \ref{sect:bincomb}. For each statistic, we first study the moments, then apply the theorem stated below, and finally address the randomness of the martingale difference sequences. 

Let us start with the statement of a Berry-Esseen type limit theorem, which can be found in \cite{Oz}. Suppose $\{X_i\}_{i\geq 1}$ is a martingale difference sequence. Let us denote $\mathbf{E}(X_i^2)$ by $\sigma_i^2$ and observe that $s_n^2= \mathbf{E}(S_n^2)=\sum_{i=1}^n \sigma_i^2$, which follows from the fact that $\mathbf{E}(X_i X_j)=0$ for $i\neq j$ as $X_i$ and $X_j$ are martingale differences. Our first assumption below is rather a technical condition. It imposes the variance of martingale differences to grow polynomially.
\begin{equation}\label{maincond}
1 \leq \liminf_n \sqrt{n}\frac{\sigma_{n+1}}{s_n}\leq\limsup_n \sqrt{n}\frac{\sigma_{n+1}}{s_n} < \infty
\end{equation}
Then the theorem is as follows.
\begin{theorem}\label{convthm}
Let $\{\mathcal{F}_{n}\}_{n\geq 1}$ be an increasing sequence of $\sigma$-fields in $\mathcal{F}$ for a probability space $(\Omega, \mathcal{F},\mathbf{P}).$ Suppose that $S_n$ is a sum of martingale differences $X_{i}$ that satisfy \eqref{maincond} and $Y_i=X_i/\sigma_i.$ If 
\begin{gather}
 \sup_i \sqrt{i}\,  \|\mathbf{E}(Y_{i}^2 | \mathcal{F}_{i-1})-1 \|_{p} < \infty  \label{bound2} \\
\sup_i \sqrt[2p']{i} \, \|\mathbf{E}(Y_{i}^3 | \mathcal{F}_{i-1}) \|_{p'} < \infty , \label{bound3} \\
\sup_i \|\mathbf{E}(Y_{i}^4 | \mathcal{F}_{i-1}) \|_\infty < \infty, \label{bound4}
\end{gather}
for some $p,p' > 1,$ then 
\begin{equation*}
\sup_{t \in \mathbb{R}} |P(S_n/s_n \leq t) - \Phi(t) | \leq \frac{C}{\sqrt{n}}
\end{equation*}
where $C$ is a constant independent of $n$ and $\Phi$ is the standard normal distribution.
\end{theorem}
 
The idea is to show that the theorem applies for any decomposition determined by the outcomes of the particle process defined in Section \ref{sect3}. We treat each statistic separately as they have distinctive features.

\subsection{The number of descents in random involutions}

Let us take $Z_n = n \left( I_n - \frac{n-1}{2} \right),$ the zero-mean stochastic process for the number of descents in random involutions. Define the central random variable $W_{i}=I_i-\frac{i-1}{2}.$  We use the method described in Section \ref{sect3} to write $Z_n$ as a random sum of martingale differences. Taking $i\geq 2,$ we denote by $X_{i,1}$ the martingale difference of the one-jump at the $i$th stage, which is $Z_{i}-Z_{i-1}$ conditioned on $\mathcal{F}_{i-1}.$ From \eqref{invjump}, we have
\begin{equation} \label{invofirdiff}
X_{i,1}=iI_i-(i-1)I_{i-1}-(i-1)=\begin{aligned}
   &  
   \begin{cases}
     W_{i-1}-\frac{i}{2},&  \text{with prob. } \frac{1}{2}+\frac{W_{i - 1}}{i},\\
       W_{i-1}+\frac{i}{2},& \text{with prob. } \frac{1}{2}-\frac{W_{i-1}}{i}.
    \end{cases}
\end{aligned}
\end{equation}
On the other hand, the two-jump at the $i$th stage is denoted by $X_{i,2}$ and given by the difference $Z_{i}-Z_{i-2}$ conditioned on $\mathcal{F}_{i-2}$ as below. 
\begin{equation}\label{invosecdiff}
X_{i,2}=iI_i - (i-2)I_{i-2} -2i +3
=2\begin{aligned}
   &   
   \begin{cases}
    W_{i-2}-\frac{i}{2},&  \text{with prob. }\frac{\left(W_{i-2}+\frac{i-1}{2}\right)^2+i-2}{i(i-1)},\\[5pt]
    W_{i-2},&  \text{with prob. } \frac{2\left(W_{i-2}+\frac{i-1}{2}\right)\left(\frac{i-1}{2}-W_{i-2}\right)-i+3}{i(i-1)},\\[5pt]
     W_{i-2}+\frac{i}{2},&  \text{with prob. } \frac{\left(\frac{i-1}{2}-W_{i-2}\right)^2+i-2}{i(i-1)}.
    \end{cases}
\end{aligned}
\end{equation}
We observe that both $\mathbf{E}[X_{i,1}|\mathcal{F}_{i-1}]$ and $\mathbf{E}[X_{i,2}|\mathcal{F}_{i-2}]$ are equal to zero. Therefore, any additive combination of them forms a martingale difference sequence.

\subsubsection{Variance} We use the generating function \eqref{invosecond} to bound the order of the variance by the falling factorial moment formula below, see Section III.2.1 of \cite{FS}. 
\begin{equation}\label{fallfact}
\mathbf{E}(X(X-1)\cdots (X-r+1)) = \frac{I_n^{(r)}(t)}{I_n(t)}\bigg|_{t=1}
\end{equation}
where $I_n^{(r)}(t)$ is the $r$th derivative of $I_n(t)$ given in  \eqref{invosecond}. We then have
\begin{align}\label{firdeg}
\begin{split}
nI_{n}'(t)
&=(n-1)I_n(t)+(t-t^2)I_{n-1}''(t)+[2+(n-3)t]I_{n-1}'(t)\\[5pt]
&\quad+t^2(1-t)^2I_{n-2}'''(t)+t(1-t)[5+(2n-9)t]I_{n-2}''(t)\\[5pt]
&\quad +[3-6t+(2n-5)(2t-3t^2)+(n-1)(1+t+(n-2)t^2)]I_{n-2}'(t)\\[5pt]
&\quad+(n-1)[1+2(n-2)t]I_{n-2}(t),
\end{split}
\end{align}
and the second derivative is
\begin{align}\label{secdeg}
\begin{split} 
nI_{n}''(t)
&=(t-t^2)I_{n-1}'''(t)+[3+(n-5)t]I_{n-1}''(t)+2(n-2)I_{n-1}'(t)
\\[5pt]
&\quad+t^2(1-t)^2I_{n-2}^{(4)}(t)+t(1-t)[7+(2n-13)t]I_{n-2}'''(t) \\[5pt]
&\quad+ [8-24t+12t^2+2(2n-5)(2t-3t^2)+(n-1)(1+t+(n-2)t^2)]I_{n-2}''(t) \\[5pt]
&\quad+[-6+(2n-5)(2-6t)+2(n-1)(1+2(n-2)t)]I_{n-2}'(t)\\[5pt]
&\quad+2(n-1)(n-2)I_{n-2}(t). 
\end{split}
\end{align}
Plugging in $t=1$ in \eqref{firdeg}, 
\begin{equation}\label{invofirstderi}
nI_{n}'(1)=(n-1)I_{n-1}'(1)+(n-1)I_n(1)+(n-2)(n-1)I_{n-2}'(1)+(n-1)(2n-3)I_{n-2}(1).
\end{equation}
We take $q_n=\frac{(n-1)i_{n-2}}{i_n},$ the probability of a two-jump at the $n$th stage. From \eqref{invofirstderi}, we have 
\begin{equation}
\mu_n=\left[\frac{n-1}{n}\mu_{n-1}+\frac{n-1}{n}\right](1-q_{n})+ \left[\frac{n-2}{n}\mu_{n-2}+\frac{2n-3}{n} \right]q_{n}.
\end{equation}
Using the formula $\mu_k = \frac{I_{k}'(1)}{I_k(1)}$ for $k=n-1$ and $n-2$, it can be verified from the above equation that $\mu_n=\frac{n-1}{2}.$ Then the second equation \eqref{secdeg} gives
\begin{align}
nI_{n}''(1)
&=(n-2)I_{n-1}''(1)+2(n-2)I_n'(t)+(n-2)(n-3)I_{n-2}''(1)
+\nonumber\\[5pt]
&\quad +2(2n-5)(n-2)I_{n-2}'(1)+2(n-1)(n-2)I_{n-2}(1).
\end{align}
Using the moment formula again, we have 
\begin{equation*}
\lambda_n= \left[\frac{n-2}{n}\lambda_{n-1}+\frac{(n-2)^2}{n}\right](1-q_{n})+ \left[\frac{(n-2)(n-3)}{n(n-1)}\lambda_{n-2}+\frac{2n^2-11n+20}{n}\right]q_n 
\end{equation*}
where $\lambda_n:=\mathbf{E}[I_n(I_n-1)].$ Recall that $I_1=0.$ Since the expression for the second moment is a second-order recurrence relation, it is difficult to obtain the exact solution. Nonetheless, we can view it as the expected value of the following stochastic process 
\begin{equation}\label{lambda}
\Lambda_n=\begin{aligned}
   &   
   \begin{cases}
   \frac{(n-2)(n-3)}{n(n-1)}\Lambda_{n-2}+\frac{2n^2-11n+20}{n},&  \text{with prob. } q_n,
   \\[5pt]
    \frac{n-2}{n}\Lambda_{n-1} + \frac{(n-2)^2}{n},&  \text{with prob. }1-q_n,
    \end{cases}
\end{aligned}
\end{equation}
with $\Lambda_1=0.$ The expectation $\mathbf{E}(\Lambda_n)=\lambda_n$ can inductively be shown from \eqref{lambda} by assuming $\mathbf{E} (\Lambda_i)=\lambda_i.$ Next, we define another process below to bound $\mathbf{E}(\Lambda_n)=\lambda_n.$ Let
\begin{equation} \label{anotherproc}
\widetilde{\Lambda}_n=\frac{n-2}{n}\widetilde{\Lambda}_{n-1} +n
\end{equation}
with $\widetilde{\Lambda}_1=\widetilde{\Lambda}_2=0.$ In order to compare it to the initial process, we write the latter as
\begin{equation*}
\widetilde{\Lambda}_n=\begin{aligned}
   &   
   \begin{cases}
   \frac{(n-2)(n-3)}{n(n-1)} \widetilde{\Lambda}_{n-2} +\frac{2n^2-3n+2}{n},&  \text{with prob. } q_n,
   \\[5pt]
    \frac{n-2}{n} \widetilde{\Lambda}_{n-1} + n,&  \text{with prob. }1-q_n,
    \end{cases}
\end{aligned}
\end{equation*}
Since both the additive terms are larger compared to \eqref{lambda} for $n\geq 3,$ the expected value of $\widetilde{\Lambda}_n$ is larger. 

Defining $\widetilde{\lambda}_n=\mathbf{E}(\widetilde{\Lambda}_n),$ we have 
\begin{equation}\label{firstinh}
\widetilde{\lambda}_n=\frac{n-2}{n} \widetilde{\lambda}_{n-1} + n
\end{equation}
from \eqref{anotherproc}. Therefore, we can bound $\lambda_n$ using the first order inhomogenous recursive relation \eqref{firstinh}. For a recursive sum obtained from the form $A_n= a_n A_{n-1} +b_n,$ we have the formula
\begin{equation}\label{inhom}
A_n= \left( \prod_{k=1}^n a_k \right) \left(A_0+\sum_{i=1}^{n}\frac{b_i}{\prod_{j=1}^i a_j}\right)
\end{equation}
in Section 2.2 of \cite{GKP89}. Let us take $a_n=\frac{n}{n+2}, b_n=n+2$ for $n \geq 1,$ and set $A_0=a_0=b_0=0.$ Then we have,
\begin{equation*}
\widetilde{\Lambda}_n = \left(\prod_{k=1}^{n-2}\frac{k}{k+2}\right) \sum_{i=1}^{n-2} \frac{j+2}{\frac{1}{(j+1)(j+2)}}=\frac{1}{n(n-1)}\sum_{j=0}^{n-1}j^3-j^2=\frac{(3n+2)(n+1)}{12}-\frac{1}{n(n-1)}.
\end{equation*}
Therefore,
\begin{equation}\label{invvar}
\textnormal{Var}(I_n) \leq \frac{(3n+2)(n+1)}{12} +\frac{n-1}{2} - \frac{(n-1)^2}{4} =\frac{17n-4}{12}=\mathcal{O}\left(n\right).
\end{equation}

\subsubsection{A bound for the fourth moment} \label{invfour}

In order to apply Theorem \ref{convthm}, we need estimates on the moments up to the fourth degree. Let us first write down the higher conditional moments of the martingale differences using \eqref{invofirdiff} and \eqref{invosecdiff}. For $i\geq 2,$ the second moments of the martingale differences are 
\begin{align} \label{invosecmom}
\begin{split}
\mathbf{E}[X_{i,1}^2 | \mathcal{F}_{i-1}] &=  \frac{i^2}{4}-W_{i-1}^2, \\
\mathbf{E}[X_{i,2}^2 | \mathcal{F}_{i-1}] &=  \frac{i(i-1)}{2}+\frac{2i(i-2)}{i-1}-\frac{2(i-2)}{i-1}W_{i-2}^2.
\end{split}
\end{align}
The third moments of the martingale differences are calculated to be
\begin{align}\label{invothimom}
\begin{split}
\mathbf{E}[X_{i,1}^3 | \mathcal{F}_{i-1}] &=  \frac{i^2}{2}W_{i-1}-2W_{i-1}^3, \\
\mathbf{E}[X_{i,2}^3 | \mathcal{F}_{i-1}] &=  \left(\frac{3}{2(i-1)}-4\right)W_{i-2}^3+ \left(i^2 + 9i - \frac{12i}{i-1}\right)W_{i-2},
\end{split}
\end{align}
and the fourth conditional moments are
\begin{align}\label{invofoumom}
\begin{split}
\mathbf{E}[X_{i,1}^4 | \mathcal{F}_{i-1}] &=  \frac{i^4}{16}+\frac{i^2}{2}W_{i-1}^2 -3 W_{i-1}^4, \\
\mathbf{E}[X_{i,2}^4 | \mathcal{F}_{i-1}] &=  \frac{48}{i-1}W_{i-2}^4-2\left(\frac{i^3+16i^2+6i-48}{i-1}\right)W_{i-2}^2+i^3(i-1)+\frac{4i^3(i-2)}{i-1}.
\end{split}
\end{align}
Next, we use Rosenthal's inequality, see Section 2.1 of \cite{HH}, to bound the fourth moment of the sum of the martingale differences. The inequality is as follows.
\begin{equation}\label{rosenthal}
 \mathbf{E}\left(\left|\sum_{i=1}^n X_i \right|^p\right) \leq C_p\left(\mathbf{E}\left(\sum_{i=1}^n \mathbf{E}[X_i^2| \mathcal{F}_{i-1}]\right)^{p/2} + \sum_{i=1}^n \mathbf{E}(|X_i|^p)\right) 
\end{equation}
where $C_p$ depends only on $p$ and  $\{X_i\}_{i\geq 1}$ is a martingale difference sequence. Relying on the recursive structure of the martingale differences in our case, which is not a property of martingales differences as such, we will bound the fourth moment of $W_n$ recursively. For the third moment condition in Theorem \ref{convthm}, we will use Lyapunov's inequality in the proof of the main theorem.
Recalling that $W_n=I_n-\mathbf{E}I_n,$ we denote the sum of martingale differences for a given composition $\mathbf{a}$ as follows.
\begin{equation}\label{invopath}
S_{n,\mathbf{a}}:=\mathbf{E}[n W_n \big| \, C_{\leq 2}(n)=\mathbf{a}]=\sum_{i} X_{\mathbf{a_i}}
\end{equation}
Moreover, we note that
\begin{equation*}
n^4 \mathbf{E}W_n^4 =  \sum_{\mathbf{a} \in \mathcal{C}_{\leq 2}(n)}\mathbf{E} S_{n,\mathbf{a}}^4 \mathbf{P}(C_{\leq 2}(n) =\mathbf{a})
\end{equation*}
by the law of total expectation. Taking $p=4$ in \eqref{rosenthal}, we have
\begin{align*}
\mathbf{E}S_{n,\mathbf{a}}^4 \leq& C_4 \left(\mathbf{E}\left(\sum_{i} \mathbf{E}[X_{\mathbf{a_i}}^2|\mathcal{F}_{i-1}]\right)^{2} + \sum_{i=1}^n \mathbf{E}X_{\mathbf{a_i}}^4\right).
\end{align*}

Therefore,
\begin{equation*}
n^4\mathbf{E}W_n^4 \leq  C_4 \max_{\mathbf{a} \in C_{\leq 2}(n)} \left( \mathbf{E}\left(\sum_{i} \mathbf{E}[X_{\mathbf{a_i}}^2|\mathcal{F}_{i-1}]\right)^{2} + \sum_{i=1}^n \mathbf{E}X_{\mathbf{a_i}}^4\right).
\end{equation*}
We then take $X_{1,1}=X_{1,2}=0$ and apply \eqref{invosecmom} and \eqref{invofoumom} to have the following upper bound:
\begin{align*}
\mathbf{E}W_n^4 & \leq \frac{C_4}{n^4} \left( \mathbf{E} \left(\sum_{i=2}^n \frac{i^2}{3} + W_{i}^2 \right)^2 + \mathbf{E}\sum_{i=2}^n \frac{48}{i-1} W_i^4+ \mathbf{E}\sum_{i=2}^n \frac{i^2}{2} W_i^2 + \sum_{i=2}^n i^4 + \mathcal{O}(i^3) \right)  \\
& \leq \frac{C_5}{n^4} \left( n^6 + n^3\sum_{i=2}^n   \mathbf{E}W_{i}^2 + \mathbf{E} \sum_{i=2}^n   W_{i}^4 + \mathbf{E}\sum_{2\leq i<j}W_i^2 W_j^2 \right) \\
&\leq \frac{C_5}{n^4} \left( n^6 + n^3\sum_{i=2}^n   \mathbf{E}W_{i}^2 + \sum_{i=2}^n  \mathbf{E} W_{i}^4 + \sum_{2\leq i<j}\sqrt{\mathbf{E}W_i^4 \mathbf{E}W_j^4}\right)
\end{align*}
where the last line is by the Cauchy-Schwarz inequality and $C_5$ is a constant independent of $n$. Assume that $\mathbf{E}W_i^4$ is an increasing function of $i.$ Since $\mathbf{E }W_i^2=\textnormal{Var}(I_i)\leq 2i$  by \eqref{invvar}, we have 
 
\begin{align*}
\mathbf{E}W_n^4 \leq  C_5 n^2 \left( 1 + \frac{1}{n}+  \frac{(n+1)\mathbf{E}(W_{n}^4)}{n^5} \right),
\end{align*}
which implies
\begin{equation}\label{4thmom}
\mathbf{E}W_n^4 = \mathcal{O}(n^2).
\end{equation}
 If $W_i$ is not increasing, letting $\mathbf{E}W_{\kappa}^4 = \max_{1 \leq i \leq n} \mathbf{E} W_i^4,$ we find that $\mathbf{E}W_{\kappa}^4$ is smaller than $C\kappa^2 \leq Cn^2$ for some constant $C$ by running the same argument for $\kappa$ in place of $n$. The conclusion \eqref{4thmom} holds true.

\subsubsection{The asymptotic normality of the conditional distribution of $I_n$}
We verify the conditions of Theorem \ref{convthm} for $I_n$ conditioned on $\mathbf{a}$ for all $\mathbf{a}$ in $\mathcal{C}_{\leq 2}(n).$ Since, below, the martingale differences will be treated in the same way regardless of whether they are obtained from one-jump or two-jumps, we will use single subscript for the martingale differences for notational convenience. For instance, the example in \eqref{onetwoone} for $Z_{12}$ would give
\begin{equation}\label{seqzero}
(X_1,\dots, X_{7})=(X_{2,2},X_{3,1},X_{5,2},X_{7,2},X_{8,1},X_{10,2},X_{12,2}).
 \end{equation}
 
Recall the definition $\sigma_i^2=\mathbf{E}(X_{i}^2).$ Since  $\mathbf{E}(X_{i}^2)$ is of order $i^2,$ the condition $\eqref{maincond}$ is satisfied. Then from the first case in \eqref{invosecmom}, we have 
\begin{equation*}
\mathbf{E}[X_{i}^2|\mathcal{F}_{i-1}]-\sigma_{i}^2= -W_{i-1}^2 + \mathbf{E}(W_{i-1}^2).
\end{equation*}
Thus,
 \begin{equation*}
 \|\mathbf{E}[X_{i}^2|\mathcal{F}_{i-1}]-\sigma_{i}^2\|_{2}^2 = \mathbf{E}\left[\left|W_{i-1}^2-\mathbf{E}(W_{i-1}^2)\right|^{2}\right] .
 \end{equation*}
Since $\mathbf{E}(W_{i-1}^2) < 2i$ for all $i$ by \eqref{invvar},
\begin{equation*}
\|\mathbf{E}[X_{i}^2|\mathcal{F}_{i-1}]-\sigma_i^2\|_{2}^2 =  \mathcal{O}(i^2).
\end{equation*} 
Taking the square root of the above expression and dividing by $\sigma_i^2,$ it is of order less than $\sqrt{i}.$ The same line of argument with different constants applies to the second case in \eqref{invosecmom}. Therefore, \eqref{bound2} is satisfied for $p=2$.

For the next condition, we bound 
\begin{equation} \label{norm3}
\|\mathbf{E}[X_i^3|\mathcal{F}_{i-1}]\|_{p}^p \leq  \mathbf{E}\left[\left|i^2 W_{i-1}+4 W_{i-1}^3 \right|^{p}\right].
\end{equation}
We take $p=4/3$ in $\eqref{norm3},$ then it is bounded by
\begin{align*}
& \mathbf{E}\left[\left(\sqrt[3]{i^2 |W_{i-1}|} +\sqrt[3]{4} \, |W_{i-1}| \right)^{4}\right] \\
 \leq &C \, \mathbf{E}(i^{8/3}|W_{i-1}|^{4/3}+i^{2}|W_{i-1}|^{2}+i^{4/3}|W_{i-1}|^{8/3}+i^{2/3}|W_{i-1}|^{10/3}+|W_{i-1}|^{4})\\
= &\mathcal{O}(i^{10/3}),
\end{align*}
where the last line follows from Lyapunov's inequality,
\begin{equation}\label{lyapunov}
\sqrt[r]{\mathbf{E}|W_{i-1}|^r} \leq \sqrt[s]{\mathbf{E}|W_{i-1}|^s} = \mathcal{O}(\sqrt{i})
\end{equation}
for $0 <r <s$, and the bounds on the second and the fourth moment, \eqref{invvar} and \eqref{4thmom} respectively. The second case in \eqref{invothimom} is treated similarly. Therefore, $\sqrt[8/3]{i}\|\mathbf{E}[X_i^3|\mathcal{F}_{i-1}]\|_{4/3}$ is also uniformly bounded. 

For the last condition, we have
\begin{align*}
\|\mathbf{E}[X_i^4|\mathcal{F}_{i-1}]\|_{\infty}&=\left\|\frac{i^4}{2}+2i^2 W_{i-1}^2 + 3 W_{i-1}^4\right\|_{\infty} = \mathcal{O}(i^4)
\end{align*}
as $|W_{i-1}|\leq \frac{i-2}{2}.$ Then since $\sigma_i$ is of order $i,$ $\|\mathbf{E}[X_i^4|\mathcal{F}_{i-1}]\|_{\infty}$ is uniformly bounded. The two-jump case is essentially the same, thereby \eqref{bound4} is satisfied.

Thus, the asymptotic normality of $S_{n,\mathbf{a}}$ follows from Theorem \ref{convthm} with an error term of order less than or equal to $n^{-1/2}$, i.e., for all $\mathbf{a} \in \mathcal{C}_{\leq 2}(n),$
\begin{equation}\label{invpathconv}
 \left|\mathbf{P}\left(\frac{S_{n,\mathbf{a}}}{s_{n,\mathbf{a}}} \leq x\right) - \Phi(x) \right|\leq \frac{C}{\sqrt{n}} 
 \end{equation} 
where $s_{n,\mathbf{a}}$ is the standard deviation of $S_{n,\mathbf{a}}.$ In the following part of the section, we will have an explicit expression for $s_{n,\mathbf{a}}^2$. We note that $C$ is independent of $n$ but depends on the limiting values in \eqref{maincond}; we refer to the proof of Theorem \ref{convthm} in Section 7 of \cite{Oz} for the details. 
 
 \subsubsection{Proof of Theorem \ref{main} for $I_n$} \label{invmain} 

We first define $z_n^2:=\textnormal{Var}(Z_n)$ recalling that $Z_n=n(I_n-\mathbf{E}(I_n)).$  We will show that the Kolmogorov distance between $\frac{Z_n}{z_n}$ and the standard normal distribution is smaller than a constant times $n^{-1/2}.$ In other words, for any $x$ in $\mathbb{R},$ we will prove that the distance between the probability below and $\Phi(x)$ is less than $C n^{-1/2}$ for some constant $C$ independent of $x$ and $n$. By definition,
 \begin{align} \label{invocompexp}
 \begin{split}
 \mathbf{P}(Z_n \leq z_n x ) = &\sum_{\mathbf{a} \in \mathcal{C}_{\leq 2}(n)} \mathbf{P}(S_{n,\mathbf{a}} \leq z_n x ) \mathbf{P} (C_{\leq 2}(n)=\mathbf{a}) \\
 =& \sum_{\mathbf{a} \in \mathcal{C}_{\leq 2}(n)} \mathbf{P}\left(\frac{S_{n,\mathbf{a}}}{s_{n,\mathbf{a}}} \leq \frac{z_n}{s_{n,\mathbf{a}}} x \right) \mathbf{P} (C_{\leq 2}(n)=\mathbf{a})
 \end{split}
 \end{align}
We start with an estimate for $s_{n,\mathbf{a}}$ to show that it is close enough to $z_n$. Note that $z_n^2$ is the average value of $s_{n,\mathbf{a}}^2,$ i.e.,
\begin{equation}\label{invovarave}
z_n^2= \sum_{\mathbf{a} \in \mathcal{C}_{\leq 2}(n)} s_{n,\mathbf{a}}^2 \mathbf{P}(C_{\leq 2}(n) =\mathbf{a})
\end{equation}
by the law of total expectation. Since the summands of $S_{n,\mathbf{a}}$ are uncorrelated, its variance decomposes as 
\begin{equation}\label{invovardec}
s_{n,\mathbf{a}}^2= \sum \sigma_{\mathbf{a_i}}^2
\end{equation}
 where $\sigma_{\mathbf{a_i}}^2=\mathbf{E}(X_{\mathbf{a_i}}^2).$  Now, considering the coefficients of $i^2$ in both types of differences in \eqref{invosecmom} and the bound \eqref{invvar}, it follows that
\begin{equation*}
\sum_{i=1}^{\floor{n/2}} \frac{(2i)^2}{2} + \mathcal{O}(i) \leq s_{n,\mathbf{a}}^2 \leq \sum_{i=1}^{n} \frac{(i)^2}{4} + \mathcal{O}(i),
 \end{equation*} 
 which implies 
 \begin{equation*}
 s_{n,\mathbf{a}}^2=\frac{n^3}{12}+\mathcal{O}(n^2).
 \end{equation*}
 Thus, $|s_{n,\mathbf{a}}^2-z_n^2|=\mathcal{O}\left(n^2\right)$ regardless of the composition of $S_{n,\mathbf{a}}.$ Then since $s_{n,\mathbf{a}}$ and $z_n$ are of order $n\sqrt{n},$ we have
 $|s_{n,\mathbf{a}}-z_n|= \mathcal{O}(\sqrt{n}),$ or, in other words,
\begin{equation}\label{invvarratio}
 \left| \frac{z_n}{s_{n,\mathbf{a}}}\right| \leq 1 + \frac{C'}{n}
 \end{equation} 
  for some constant $C'.$ 
\begin{remark}
In fact, $S_{n,\mathbf{a}}$ typically consists of two-jumps. If we define $T_n \hspace*{-0.5mm} =\sum_{i=1}^n \mathcal{B}\left(\frac{1}{\sqrt{i}},1,0\right)$ as in Section \ref{sect:app}, which gives a binary approximation to the number of one-jumps in \eqref{invjump} and bound it from above by the estimate \eqref{invratio}, then we can show by Chebyshev's inequality that the martingale difference sequence has less than $4\sqrt{n}$ one-jumps with probability at least $1-n^{-1/2}$.
\end{remark} 

Now, it follows from \eqref{invocompexp}, \eqref{invvarratio} and \eqref{invpathconv} that
 \begin{align*}
 \begin{split}
 \mathbf{P}(Z_n \leq z_n x ) \leq &\sum_{\mathbf{a} \in \mathcal{C}_{\leq 2}(n)} \mathbf{P}\left(\frac{S_{n,\mathbf{a}}}{s_{n,\mathbf{a}}} \leq \left(1+\frac{C'}{\sqrt{n}}\right) x \right)  \\
 \leq & \left|  \Phi\left( 1+\frac{C'}{n}  \right)x \right| + Cn^{-1/2}.
\end{split} 
 \end{align*}
Finally, we use the following estimate,
 \begin{equation}\label{normalest}
 \left| \Phi\left( \left( 1+a \right)x \right) - \Phi\left( x\right) \right| \leq \frac{e^{-1/2}}{\sqrt{2 \pi}}a
 \end{equation}
for $a \in \mathbb{R},$ to conclude that
 \begin{equation*}
 \left|\mathbf{P}\left(\frac{Z_n}{z_n} \leq x \right) -\Phi(x)\right| \leq \left(\frac{e^{-1/2}C'}{\sqrt{2 \pi n } }+C\right)n^{-1/2}.
 \end{equation*}
 \bbox

\subsection{The number of descents in random derangements}

Let $Z_n=(n-1)(R_n-\mu_n)$ where $\mu_n:=  \mathbf{E}(R_n).$  We will show that $Z_n$ converges to the normal distribution as in the previous section. However, we have two major differences here. We start with noting that neither $Z_{i}-Z_{i-1}$ nor $Z_{i}-Z_{i-2}$ render a martingale difference. Although the conditional expectation for the former is deterministic, the latter cannot even be corrected by an additive term. Yet, we observe that the conditional expectation of $Z_{i}-\frac{i-2}{i-3}Z_{i-2}$ is constant. We thus need an adjustment factor in the decomposition of $Z_n$ into martingale differences, which is the first major difference. For example, if the composition is of the form $(\dots,2,1,1,1,2,2,1),$ then the differences are
\begin{align*}
Z_{n+2} =&(n+1)R_{n+2}-n R_{n+1}+ \\
&n R_{n+1} - (n-1) R_{n-1} + \\
&\frac{n-1}{n-2} \left[(n-2)R_{n-1} - (n-3) R_{n-3} \right]+ \\
&\left(\frac{n-3}{n-4}\right)\left(\frac{n-1}{n-2}\right) \left[ (n-4)R_{n-3} - (n-5) R_{n-4}\right] + \\
&\left(\frac{n-3}{n-4}\right)\left(\frac{n-1}{n-2}\right) \left[ (n-5)R_{n-4} - (n-6) R_{n-5}\right]+ \ldots
\end{align*}
Let us denote the product of adjustment factors down to the $i$th term by  
\begin{equation*}
\Gamma_{\mathbf{a}}(i)=\prod_{k>i:a_k=2} \frac{k}{k-1}.
\end{equation*}
The following bound will be useful, 
\begin{equation}\label{wallis}
\Gamma_{\mathbf{a}}(i) \leq \Gamma_{(2,\dots,2)}(1) = \prod_{k=1}^{\floor{n/2}} \frac{2k}{2k-1} \leq \sqrt{\frac{\pi n}{2}}\left(1+\mathcal{O}(n^{-1}) \right), 
\end{equation}
which is a consequence of Wallis' product formula, see 6.1.49 of \cite{A64}.

The second major difference is the deterministic additive term, denoted by $\alpha_{i,1}$ for one-jumps and $\alpha_{i,2}$ for two-jumps below. From \eqref{derjump}, we have
\begin{equation*}
Z_i-Z_{i-1}=X_{i,1}+\alpha_{i,1}
\end{equation*}
where 
\begin{equation}\label{derdiff1}
X_{i,1}= \begin{aligned}
   &   
   \begin{cases}
     R_{i-1}-i+2,&  \text{with prob. } \frac{R_{i-1}+1}{(i-1)},\\[5pt]
     R_{i-1}+1,&  \text{with prob. } \frac{i-2-R_{i-1}}{(i-1)}
    \end{cases}
\end{aligned}
\end{equation}
and $\alpha_{i,1}:=i-2-(i-1)\mu_{i}+(i-2)\mu_{i-1}.$ We split the difference in this particular way to guarantee $\mathbf{E}[X_{i,1}|\mathcal{F}_{i-1}]=0.$ Similarly,
\begin{equation} \label{derdeux}
Z_{i}-\frac{i-2}{i-3}Z_{i-2}=X_{i,2}+\alpha_{i,2}
\end{equation}
 where 
\begin{equation}\label{derdiff2}
X_{i,2}=\begin{aligned}
   &   
   \begin{cases}
     R_{i-2}-i+2,&  \text{with prob. } \frac{R_{i-2}+1}{(i-1)},\\[5pt]
     R_{i-2}+1,&  \text{with prob. } \frac{i-2-R_{i-2}}{(i-1)}
    \end{cases}
\end{aligned}
\end{equation}
and $\alpha_{i,2}:=2i-3-(i-1)\mu_{i}+(i-2)\mu_{i-2}.$ 

Thus, considering both differences, we have
\begin{equation}\label{formder}
\mathbf{E}[Z_n | C_{\leq 2}(n)=\mathbf{a}]=\sum_i \Gamma_{\mathbf{a}}(i)  \left( X_{\mathbf{a_i}}+ \alpha_{\mathbf{a_i}} \right)
\end{equation}
compared to \eqref{forminv}. Let us use the notation $S_{n, \mathbf{a}}$ and $\alpha_{n,\mathbf{a}}$ for 
$\sum_i \Gamma_{\mathbf{a}}(i)  X_{\mathbf{a_i}}$ and $ \sum_i \Gamma_{\mathbf{a}}(i) \alpha_{\mathbf{a_i}}$ respectively.

\subsubsection{Moments}  We use the factorial moment formula \eqref{fallfact} to estimate the moments by the generating function \eqref{deragen}. The first derivative of $D_n(t)$ is
\begin{align} \label{derivder}
D_n'(t)=&(-1)^n (n-1)t^{n-2}+(n-1)D_{n-1}(t)+(1+(n-1)t)D_{n-1}'(t) \\
&+(1-2t)D_{n-1}'(t) + t(1-t)D_{n-1}''(t). \notag
\end{align}
 We denote $D_n'(1)$ by $d_n'$ for notational ease. We evaluate the above expression at $t=1,$
\begin{align*}
d_n'&=(-1)^n(n-1)+(n-1)d_{n-1}+ (n-1)d_{n-1}' \\
&=(-1)^n(n-1)+(n-1)((n-1)d_{n-2}+(-1)^{n-1})+ (n-1)d_{n-1}'\\
&=(n-1)^2d_{n-2}+ (n-1)d_{n-1}'.
\end{align*}
The following recurrence relation for the expected value follows from \eqref{fallfact}.
\begin{equation*}
\mu_n=\frac{(n-1)d_{n-1}}{d_n} \mu_{n-1} + \frac{(n-1)^2 d_{n-2}}{d_n}.
\end{equation*}
By the recusive sum formula \eqref{inhom},
\begin{equation*}
\mu_n = \left( \prod_{k=1}^{n-1} \frac{k d_k}{d_{k+1}} \right) \left( \sum_{m=1}^{n-1} m^2 \frac{d_{m-1}}{d_{m+1}}  \middle/ \prod_{k=1}^{m} \frac{k d_k}{d_{k+1}} \right) 
= \frac{(n-1)!}{d_n} \sum_{m=1}^{n-1} \frac{m^2}{m!} d_{m-1}.
\end{equation*}
Then using the closed-form expression for the number of derangements \eqref{derclose}, we have
\begin{equation*}
\sum_{m=1}^{n-1} \frac{m^2}{m!} d_{m-1} = \sum_{j=0}^{n-1} \left( \sum_{k=j+1}^{n-1} k \right) \frac{(-1)^j}{j!} = \sum_{j=0}^{n-1}\left( \frac{n(n-1)}{2}-\frac{j(j+1)}{2}\right)\frac{(-1)^j}{j!}.
\end{equation*}
The contribution of the first term in the paranthesis on the right-hand side is 
\begin{equation*}
\frac{n(n-1)}{2}\sum_{j=0}^{n-1}\frac{(-1)^j}{j!}=\frac{n(n-1)}{2}\left(\sum_{j=0}^{n}\frac{(-1)^j}{j!}+\frac{(-1)^{n+1}}{n!}\right)=\frac{n(n-1)}{2}\frac{d_n}{n!} + \mathcal{O}\left(\frac{1}{(n-2)!}\right). 
\end{equation*}
considering the error term in the Taylor expansion of $e^{-x}.$ We then evaluate the second term as follows.
\begin{align*}
-\frac{1}{2} \sum_{j=0}^{n-1} j(j+1)\frac{(-1)^j}{j!} =& \frac{1}{2}\frac{d}{dx}\bigg|_{x=1} \sum_{j=0}^{n-2} \frac{(-1)^{j}}{j!}x^{j+2}\\
=&\frac{d}{dx}\bigg|_{x=1} \frac{x^2 e^{-x}}{2} -  \frac{1}{2} \sum_{j=n}^{\infty} (j+1) \frac{(-1)^{j-1}}{(j-1)!}  \\
=& \frac{e^{-1}}{2} +  \mathcal{O}\left(\frac{1}{(n-2)!}\right).
\end{align*}
The error term estimate also implies $e^{-1}/\sum_{j=0}^{n}\frac{(-1)^j}{j!}=\mathcal{O}\left(\frac{1}{(n-2)!}\right).$ Therefore, we have
\begin{equation}\label{dermom}
\mathbf{E}R_n= \frac{n-1}{2} + \frac{1}{2n} + o\left(e^{-n}\right).
\end{equation}

Regarding the second moment, we differentiate $\eqref{derivder}$ to obtain
\begin{align*}
d_n''=&(-1)^n (n-1)(n-2) + (2n-4)d_{n-1}'+ (n-2)d_{n-1}'',
\end{align*}
where $d_n'':=D_n(1)''.$ Then the recurrence relation
\begin{equation}\label{dersecfall}
\lambda_n= \frac{(n-2)d_{n-1}}{d_n} \lambda_{n-1} + \frac{(2n-4)d_{n-1} \mu_{n-1}}{d_n} +\frac{(-1)^n(n-1)(n-2)}{d_n}
\end{equation}
for $\lambda_n=\mathbf{E}(R_n(R_n-1))$ follows from \eqref{fallfact}. We apply the same formulas \eqref{inhom} and \eqref{derclose} to \eqref{dersecfall} in this case.
\begin{align*}
\lambda_n =& \left( \prod_{k=2}^{n-1} \frac{(k-1) d_k}{d_{k+1}} \right) \left( \sum_{m=2}^{n-1} \frac{(2m-2) d_{m}\mu_m+(-1)^{m+1}m(m-1)}{d_{m+1}}  \middle/ \prod_{k=2}^{m} \frac{(k-1) d_k}{d_{k+1}} \right) \\
=& \frac{(n-2)!}{d_n}\left( \sum_{m=2}^{n-1} \frac{(2m-2)d_m \mu_m}{(m-1)!} + \frac{(n-2)!}{d_n} \sum_{m=2}^{n-1} \frac{(-1)^{m+1}m(m-1)}{(m-1)!} \right)\\
=& \frac{(n-2)!}{d_n} \left(\sum_{m=2}^{n-1} \frac{(2m-2)m!}{(m-1)!}  \left( \frac{m-1}{2} + \frac{1}{2m} + o\left(e^{-m}\right)\right) \sum_{i=0}^m\frac{(-1)^i}{i!}- e^{-1} + o(e^{-n})  \right)\\
=& \frac{1}{n(n-1)} \sum_{m=2}^{n-1} m(m-1)^2 \left(\sum_{i=0}^m\frac{(-1)^i}{i!} \middle/ \sum_{i=0}^n\frac{(-1)^i}{i!}\right) + o(n)\\
=& \frac{1}{n(n-1)} \sum_{m=2}^{n-1} (m^3-2m^2+m) \left(1+\mathcal{O}\left(\frac{1}{(m-2)!}\right)\right) + o(n)\\
=& \frac{3n^2-11n}{12} + o(n).
\end{align*}
Finally, we obtain the order of the variance from the second factorial moment.
\begin{equation}\label{varder}
\textnormal{Var}(R_n)=\frac{3n^2-11n}{12}+\frac{n-1}{2}-\frac{(n-1)^2}{4}+o(n)=\frac{n}{12} +o(n).
\end{equation}

\subsubsection{The deterministic term} \label{sec:det} We consider the sum of the deterministic terms in \eqref{formder}, $\alpha_{n,\mathbf{a}},$  as a random variable over the compositions in $\mathcal{C}_{\leq 2}(n)$ and bound its moments. We will eventually show that it does not change the asymptotic distribution but possibly affects the rate of convergence in the limit. Define
\begin{equation*}
\alpha_{i}=\begin{cases}
\alpha_{i,2} &\textnormal{ with prob. } \frac{(i-1)d_{i-2}}{d_i} \\
\alpha_{i,1} &\textnormal{ with prob. } \frac{(i-1)d_{i-1}}{d_i}.
\end{cases}
\end{equation*}
From the formula \eqref{dermom} for $\mu_n$,
\begin{equation*}
\alpha_{i}=\begin{cases}
\frac{i-1}{2} + \mathcal{O}(i^{-1}) &\textnormal{ with prob. } \frac{(i-1)d_{i-2}}{d_i} \\
-\frac{1}{2} + \mathcal{O}(i^{-1}) &\textnormal{ with prob. } \frac{(i-1)d_{i-1}}{d_i}.
\end{cases}
\end{equation*}
Since $\frac{d_{i-2}}{d_{i-1}}=\frac{1}{i-1}+ o(e^{-i})$ by \eqref{derclose}, $\alpha_i$ is bounded above by $\mathcal{B}\left(\frac{1}{i}, \frac{i}{2}, -\frac{1}{2}\right)$ up to an exponentially small error term. We then bound the adjustment factor for the $i$th term in the sum \eqref{formder},
\begin{equation}\label{wallisi}
\Gamma_{\mathbf{a}}(i) \leq \prod_{i=1}^{\floor{\frac{n-k}{2}}} \frac{n-2i+1}{n-2i} \leq \sqrt{\frac{n}{i}}
\end{equation}
by \eqref{wallis}.
Let $B_i = \left|\mathcal{B}\left(\frac{1}{i}, \frac{i}{2}, -\frac{1}{2}\right)\right|$ and define
\begin{equation}
U_n=\sum_{i=1}^n \sqrt{\frac{n}{i}}B_i.
\end{equation}
Observe that 
\begin{equation}
\mathbf{E}B_i=1-\frac{1}{2i} \, \textnormal{ and } \, \mathbf{E}B_i^k= \frac{i^{k-1}+1}{2^k}+\mathcal{O}(i^{-1}) \textnormal{ for all integer }k\geq 2,
\end{equation}
which give
\begin{gather} 
\mathbf{E}U_n \leq \sqrt{n}\sum_{i=1}^n \frac{1}{\sqrt{i}}\leq 2 n +\mathcal{O}(\sqrt{n}), \label{youn} \\
\mathbf{E}U_n^2\leq n \left(\sum_{i=1}^n \frac{i+1}{4i} + \sum_{1\leq i, j\leq n} \frac{1}{\sqrt{ij}}\right) \leq 5 n^2, \label{youn2} \\
\mathbf{E}U_n^3\leq n \sqrt{n}\left(\sum_{i=1}^n \frac{i^2+1}{8i\sqrt{i}} + 3\sum_{1 \leq i,j\leq n}  \frac{i+1}{4i\sqrt{j}} +  \sum_{1\leq i,j,k\leq n} \frac{1}{\sqrt{ijk}} \right) \leq 10n^3, \notag \\
\mathbf{E}U_n^4\leq n^2 \left(\sum_{i=1}^n \frac{i^3+1}{16i^2} + 4 \hspace*{-1mm} \sum_{1\leq i,j \leq n}  \frac{i^2+1}{8i\sqrt{ij}} +  6 \hspace*{-1mm}\sum_{1\leq i,j \leq n} \frac{(i+1)(j+1)}{16ij}+ \hspace*{-2mm} \sum_{1\leq i,j,k,l\leq n} \frac{1}{\sqrt{ijkl}} \right) \leq 20n^4 \notag
\end{gather}
by expanding the powers of $U_n$ and using the independence of its summands.
Since $|\alpha_{n,\mathbf{a}}| \leq \psi(U_n) \leq U_n,$ we have
\begin{equation}\label{detmom}
 \mathbf{E}\alpha_{n,\mathbf{a}}^k = \mathcal{O}(n^k)
 \end{equation} 
for $k=1,2,3$ or $4.$

\subsubsection{A bound for the fourth moment} 
Let us define $W_{i,1}=R_{i-1}-\frac{i-3}{2}$ and $W_{i,2}=R_{i-2}-\frac{i-3}{2}$ and also the central random variable $\overline{W}_i=R_i - \mathbf{E}(R_i).$ It follows from \eqref{dermom} that 
\begin{align}\label{derapp}
\begin{split}
W_{i,1} &= \overline{W}_i-\frac{1}{2}+\frac{1}{2(i-1)}+\mathcal{O}(e^{-i}), \\
W_{i,2} &= \overline{W}_i+\frac{1}{2(i-2)}+\mathcal{O}(e^{-i}).
\end{split}
\end{align}
By \eqref{derdiff1} and \eqref{derdiff2}, the second moments of the martingale differences are 
\begin{equation} \label{dersecmom}
\mathbf{E}[X_{i,*}^2 | \mathcal{F}_{i-1}] =  \frac{(i-1)^2}{4}-W_{i,*}^2, 
\end{equation}
for $W_{i,*}=W_{i,1}$ or $W_{i,2}.$ Similarly, the third moments of the martingale differences are 
\begin{equation}\label{derthimom}
\mathbf{E}[X_{i,*}^3 | \mathcal{F}_{i-1}] =  \frac{(i-1)^2}{2}W_{i,*}-2W_{i,*}^3
\end{equation}
and the fourth moments are
\begin{equation}\label{derfoumom}
\mathbf{E}[X_{i,*}^4 | \mathcal{F}_{i-1}] =  \frac{(i-1)^4}{16}+\frac{(i-1)^2}{2}W_{i,*}^2 -3 W_{i,*}^4.
\end{equation}
The sum of martingale differences takes the form 
\begin{equation}\label{derpath}
S_{n,\mathbf{a}}=\sum_i \Gamma_{\mathbf{a}}(i)  X_{\mathbf{a_i}} =(n-1)\overline{W}_n-\alpha_{n,\mathbf{a}},
\end{equation}
which is obtained from \eqref{formder}. Let us first show that  $\alpha_{n, \mathbf{a}}$ is insignificant regarding the order of the bound on the fourth moment. Denote the expectation with respect to compositions in $\mathcal{C}_{\leq 2}(n)$ by $\mathbf{E}_{\mathbf{a}}.$ Then
\begin{align}\label{wbar}
\begin{split}
(n-1)^4\mathbf{E}\overline{W}_n^4 &= \sum_{\mathbf{a} \in \mathcal{C}_{\leq 2}(n)}\mathbf{E} \left(S_{n,\mathbf{a}}+\alpha_{n,\mathbf{a}}\right)^4 \mathbf{P}(C_{\leq 2}(n) =\mathbf{a}) \\
&= \mathbf{E}_{\mathbf{a}}\, \mathbf{E} \left(S_{n,\mathbf{a}}+\alpha_{n,\mathbf{a}}\right)^4\\
&= \mathbf{E} \,\mathbf{E}_{\mathbf{a}}  \left(S_{n,\mathbf{a}}+\alpha_{n,\mathbf{a}}\right)^4 \\
& \leq \mathbf{E}  \left(\sqrt[4]{\mathbf{E}_{\mathbf{a}}(S_{n,\mathbf{a}}^4)} + \sqrt[4]{\mathbf{E}_{\mathbf{a}}\left(\alpha_{n,\mathbf{a}}^4\right) }\right)^4 \\
& = \mathbf{E}  \left(\sqrt[4]{\mathbf{E}_{\mathbf{a}}(S_{n,\mathbf{a}}^4)} + \mathcal{O}(n) \right)^4
\end{split}
\end{align}
where we  used Minkowski's inequality in the fourth line and \eqref{detmom} at the bottom. 

Next, we bound the fourth moment of the sum of martingale differences. Taking $p=4$ in \eqref{rosenthal}, we have
\begin{align*}
\mathbf{E}S_{n,\mathbf{a}}^4  \leq& C_4\left( \mathbf{E}\left(\sum_{i}  \Gamma_{\mathbf{a}}(i)^2 \mathbf{E}[X_{\mathbf{a_i}}^2|\mathcal{F}_{i-1}]\right)^{2} + \sum_{i=1}^n \Gamma_{\mathbf{a}}(i)^4 \mathbf{E}X_{\mathbf{a_i}}^4\right)
\end{align*}
for any composition $\mathbf{a} \in C_{\leq 2}(n).$ By the upper bound on the adjustment factor \eqref{wallisi},  
\begin{align*}
\mathbf{E}S_{n,\mathbf{a}}^4 \leq \max_{\mathbf{a} \in \mathcal{C}_{\leq 2}(n)} C_4 \left( \mathbf{E}\left(\sum_{i} \frac{n}{i} X_{\mathbf{a_i}}^2\right)^{2} + \sum_{i=1}^n \frac{n^2}{i^2} \mathbf{E}X_{\mathbf{a_i}}^4 \right).
\end{align*}
Then the conditional moments in \eqref{dersecmom} and \eqref{derfoumom} give
\begin{align*}
\mathbf{E}S_{n,\mathbf{a}}^4 \leq & C_4 n^2 \left(  \mathbf{E} \left(\sum_{i=1}^n \frac{(i-1)^2}{4i} + \frac{W_{i,*}^2}{i} \right)^2 +\sum_{i=1}^n \frac{(i-1)^4}{16i^2} +\frac{(i-1)^2\mathbf{E}W_{i,*}^2}{2i^2}  \right) \\
\leq & 2  C_4 n^2 \left(\frac{n^4}{64}+ \frac{n^2}{8}\mathbf{E} \sum_{i=1}^n \frac{W_{i,*}^2}{i} + \mathbf{E} \sum_{i=1}^n   \frac{W_{i,*}^4}{i^2} + 2\mathbf{E}\sum_{i<j}\frac{W_{i,*}^2 W_{j,*}^2}{ij} \right) \\
\leq & 2 C_4 n^2 \left(\frac{n^4}{64}+ \frac{n^2}{8} \sum_{i=1}^n \frac{\mathbf{E}W_{i,*}^2}{i} +  \sum_{i=1}^n   \frac{\mathbf{E}W_{i,*}^4}{i^2} + 2\sum_{i<j}\frac{\sqrt{\mathbf{E}W_{i,*}^4 \mathbf{E}W_{j,*}^4}}{ij} \right).
\end{align*}
We observe that $\mathbf{E}(W_{i,*}^2)$ is less than $2i$ from the definition of $W_{i,*}$ and \eqref{varder}. We then make the assumption that $\mathbf{E}(W_{i,*}^4)$ is an increasing function of $i;$ otherwise, the argument in the end of Section \ref{invfour} applies in the same manner. Therefore,
\begin{align}\label{4thmomend}
\begin{split}
\mathbf{E}S_{n,\mathbf{a}}^4 \leq & 2 C_4 n^2 \left(\frac{n^4}{24^2}+2n^3+ \mathbf{E}W_{n,*}^4 \sum_{i=1}^n  \frac{1}{i^2} + 2 \mathbf{E}W_{n,*}^4 \sum_{i<j}\frac{1}{ij} \right)  \\
\leq & 2 C_4 n^2  \left(\frac{n^4}{24}+ \mathbf{E}W_{n,*}^4 \left( \log^2 n + \frac{1}{n}\right)\right) 
\end{split}
\end{align}
Since $\left| W_{n,*}- \overline{W}_n \right| \leq \frac{1}{2}$ by \eqref{derapp} and the second moments of both $W_{i,*}$ and $\overline{W}_i$ are larger than of order of $n,$ we have 
\begin{equation*}
\frac{\mathbf{E}W_{n,*}^4}{\mathbf{E}\overline{W}_n^4}= \mathcal{O}(1).
\end{equation*}
We can thus combine \eqref{wbar} and \eqref{4thmomend} by interchanging $W_{i,*}$ and $\overline{W}_i$ to arrive at
\begin{equation}\label{4thmomd}
\mathbf{E}W_{n,*}^4 = \mathcal{O}(n^2).
\end{equation}

\subsubsection{The asymptotic normality of the conditional distribution of $R_n$}

We verify the conditions of Theorem \ref{convthm} for $R_n$ conditioned on $\mathbf{a}$ for all $\mathbf{a}$ in $\mathcal{C}_{\leq 2}(n).$ We use the single subscript notation in the same way as in \eqref{seqzero}, so that we will write $X_i$ instead of $X_{i,*},$ and $W_i$ instead of  $W_{i,*}.$

We first verify $\eqref{maincond}.$ Since $\Gamma_{\mathbf{a}}(i)$ is deterministic given $\mathbf{a}$ in $\mathcal{C}_{\leq 2}(n),$ $\sigma_i^2=\Gamma_{\mathbf{a}}(i)^2\mathbf{E}[X_{i}^2].$ Therefore, $\sigma_i^2$ is of order at most $ni$ by \eqref{wallisi}. Then $s_{n,\mathbf{a}}^2=\sum_i \sigma_i^2$ is of order $n^3$ at most, which verifies the condition. We also observe that $\Gamma_{\mathbf{a}}(i)$ being deterministic, its contribution is cancelled when $X_i$ is divided by $\sigma_i$ to obtain $Y_i$ in the statement of Theorem \ref{convthm}. Thus, we omit $\Gamma_{\mathbf{a}}(i)$ in the rest of the proof.

For the second moment condition in the theorem, we consider \eqref{dersecmom} to have
 \begin{equation*}
 \|\mathbf{E}[X_{i}^2|\mathcal{F}_{i-1}]-\sigma_{i}^2\|_{2}^2 = \mathbf{E}\left[\left|\mathbf{E}(W_{i-1}^2) - W_{i-1}^2\right|^{2}\right]=\mathcal{O}(i^2),
 \end{equation*}
which follows from \eqref{varder} and \eqref{4thmomd}. Taking the square root of the above expression and then dividing it by $\sigma_i^2,$ we see that it is of order less than $\sqrt{i}.$ Therefore, \eqref{bound2} is satisfied for $p=2$.

For the third moment condition, we bound 
\begin{equation} \label{dernorm3}
\|\mathbf{E}[X_i^3|\mathcal{F}_{i-1}]\|_{p}^p \leq  \mathbf{E}\left[\left|\frac{i^2}{4} W_{i-1}+2 W_{i-1}^3 \right|^{p}\right]
\end{equation}
from \eqref{derthimom}. Take $p=4/3$, then the expression is bounded by
\begin{align*}
& \mathbf{E}\left[\left(\sqrt[3]{\frac{1}{4}|W_{i-1}|} +\sqrt[3]{2} \, |W_{i-1}| \right)^{4}\right] \\
 \leq &C \, \mathbf{E}(i^{8/3}|W_{i-1}|^{4/3}+i^{2}|W_{i-1}|^{2}+i^{4/3}|W_{i-1}|^{8/3}+i^{2/3}|W_{i-1}|^{10/3}+|W_{i-1}|^{4})\\
= &\mathcal{O}(i^{10/3}),
\end{align*}
where the last line follows from Lyapunov's inequality \eqref{lyapunov}, \eqref{varder} and \eqref{4thmomd}. Therefore, $\sqrt[8/3]{i}\|\mathbf{E}[Y_i^3|\mathcal{F}_{i-1}]\|_{4/3}$ is also uniformly bounded. 

The last condition involves \eqref{derfoumom}, for which we have
\begin{equation*}
\|\mathbf{E}[X_i^4|\mathcal{F}_{i-1}]\|_{\infty} \leq \left\|\frac{i^4}{16}+\frac{i^2}{2} W_{i-1}^2 + 3 W_{i-1}^4\right\|_{\infty} = \mathcal{O}(i^4)
\end{equation*}
as $|W_{i-1}|\leq \frac{i-2}{2}.$ Then since $\sigma_i$ is of order $i,$ $\|\mathbf{E}[Y_i^4|\mathcal{F}_{i-1}]\|_{\infty}$ is uniformly bounded, \eqref{bound4} is also satisfied.

Therefore, by Theorem \ref{convthm}, we have for all $\mathbf{a} \in \mathcal{C}_{\leq 2}(n),$
\begin{equation}\label{derpathconv}
 \left|\mathbf{P}\left(\frac{S_{n,\mathbf{a}}}{s_{n,\mathbf{a}}} \leq x\right) - \Phi(x) \right|\leq \frac{C}{\sqrt{n}} 
 \end{equation} 
where $s_{n,\mathbf{a}}$ is the standard deviation of $S_{n,\mathbf{a}}$ and $C$ is a constant independent of $n.$

\subsubsection{Proof of Theorem \ref{main} for $R_n$} \label{dermain} 

 Recall that  $Z_n=(n-1)(R_n-\mathbf{E}(I_n))$ and let $z_n:=\textnormal{Var}(Z_n).$ For any $x$ in $\mathbb{R},$ 
 \begin{align} \label{dercompexp}
 \begin{split}
 \mathbf{P}(Z_n \leq z_n x ) = &\sum_{\mathbf{a} \in \mathcal{C}_{\leq 2}(n)} \mathbf{P}(S_{n,\mathbf{a}} + \alpha_{n,\mathbf{a}} \leq z_n x ) \mathbf{P} (C_{\leq 2}(n)=\mathbf{a}) \\
 =& \sum_{\mathbf{a} \in \mathcal{C}_{\leq 2}(n)} \mathbf{P}\left(\frac{S_{n,\mathbf{a}}+\alpha_{n,\mathbf{a}}}{s_{n,\mathbf{a}}} \leq \frac{z_n}{s_{n,\mathbf{a}}} x \right) \mathbf{P} (C_{\leq 2}(n)=\mathbf{a}) \\
 = & \sum_{\mathbf{a} \in \mathcal{C}_{\leq 2}(n)} \Phi\left( \frac{z_n}{s_{n,\mathbf{a}}} x -\frac{\alpha_{n,\mathbf{a}}}{s_{n,\mathbf{a}}} \right) \mathbf{P} (C_{\leq 2}(n)=\mathbf{a}) +\mathcal{O}(n^{-1/2})
 \end{split}
 \end{align}
 by \eqref{derpathconv}.

 Define $T_n=\sum_{i=1}^n \mathcal{B}\left(\frac{i-1}{i},0,1\right),$ which gives an upper bound on the number of two-jumps in \eqref{derjump}. We have 
\begin{equation*}
\mathbf{E}T_n=\sum_{i=1}^n \frac{1}{i} = \log{n} + \mathcal{O}(1) \textnormal{ and } \textnormal{ Var}(T_n) = \sum_{i=1}^n \frac{i-1}{i^2} =  \log n + \mathcal{O}(1).
\end{equation*}
By  Chebyshev's inequality,
\begin{equation*}
\mathbf{P}\left(\left| T_n - \log n \right| > \sqrt[4]{n \log^2 n}\right) \leq \frac{1}{\sqrt{n}}.
\end{equation*}
Since $\psi(T_n) \leq T_n$ by definition, $\psi(T_n)$ is less than $\sqrt[4]{n} \log n$ with probability at least $1- n^{-1/2}.$ Using the variance decomposition \eqref{invovardec}, we have 

\begin{equation*}
\sum_{i=1}^{\floor{n-\sqrt[4]{n} \log n}} \frac{i^2}{4} + \mathcal{O}(i) \leq s_{n,\mathbf{a}}^2 \leq \sum_{i=1}^{n} \frac{i^2}{4} + \mathcal{O}(i) \quad \textnormal{  with prob. at least }1-n^{-1/2}
 \end{equation*} 
by \eqref{dersecmom} and \eqref{varder}.

Therefore,
\begin{equation} \label{dersn}
 s_{n,\mathbf{a}}^2=\frac{n^3}{12}+\mathcal{O}(n^2) \quad \textnormal{  with prob. at least }1-n^{-1/2},
 \end{equation} 
 otherwise it is bounded above by $\frac{n^3}{8}$ in the case $\mathbf{a}=(2,2,\dots,2).$ Thus,
 \begin{equation*}
\frac{n^3}{12}+ \mathcal{O}(n^2)  \leq z_n^2 \leq (1-n^{-1/2}) \frac{n^3}{12} + n^{-1/2}\frac{n^3}{8}= \frac{n^3}{12}+ \mathcal{O}(n^2\sqrt{n})
 \end{equation*}
 by \eqref {invovarave}. This eventually shows that  $\left| \frac{z_n}{s_{n,\mathbf{a}}}\right| \leq 1 + \frac{C'}{\sqrt{n}}$ for some constant $C'.$ Then the bound in \eqref{dercompexp} becomes
 \begin{align}\label{secbound}
 \begin{split}
 \mathbf{P}(Z_n \leq z_n x ) \leq & \sum_{\mathbf{a} \in \mathcal{C}_{\leq 2}(n)} \Phi\left( \left( 1+\frac{C'}{\sqrt{n}} \right)  x -\frac{\alpha_{n,\mathbf{a}}}{s_{n,\mathbf{a}}} \right) \mathbf{P} (C_{\leq 2}(n)=\mathbf{a}) + \mathcal{O}(n^{-1/2})
\end{split} 
 \end{align}
We then turn to the deterministic term. Another application of Chebyshev's inequality gives
  \begin{equation}\label{minmax1}
\mathbf{P}\left(\left| U_n - 2n  \right| > \sqrt{5}n^{1+k}\right) \leq n^{-2k}
\end{equation}
by \eqref{youn} and \eqref{youn2}. As noted in Section \ref{sec:det}, $|\alpha_{n,\mathbf{a}}|\leq U_n,$ therefore $|\alpha_{n,\mathbf{a}}|\leq \sqrt{5}n^{1+k} + \mathcal{O}(n).$ So,
\begin{equation}\label{minmax2}
\left|\frac{\alpha_{n,\mathbf{a}}}{s_{n,\mathbf{a}}}\right|=2 \sqrt{15}n^{k-1/2} + \mathcal{O}\left(n^{k-3/4}\right)
\end{equation}
by \eqref{dersn}.
Considering the estimate 
\begin{equation*}
\left|\Phi\left(x+ a\right)-\Phi\left(x\right)\right|\leq \frac{a}{\sqrt{2 \pi}}
\end{equation*}
for all $x \in \mathbb{R}$ in \eqref{secbound}, \eqref{minmax2} becomes a part of the error term. We then optimize the exponent that appears in both  \eqref{minmax1} and \eqref{minmax2}. We have $\min_k \left(\max \left\{k-\frac{1}{2},-2k\right\}\right)=-\frac{1}{3}.$ Therefore, we can further improve \eqref{secbound} to have
 \begin{align*}
 \begin{split}
\left| \mathbf{P}(Z_n \leq z_n x ) -  \Phi\left( \left( 1+\frac{C'}{\sqrt{n}}  \right)x\right) \right| \leq 2\sqrt{15}n^{-1/3}.
\end{split} 
 \end{align*}
Finally, 
 \begin{equation*}
 \left|\mathbf{P}\left(\frac{Z_n}{z_n} \leq x \right) -\Phi(x)\right| \leq C n^{-1/3}
 \end{equation*}
 by \eqref{normalest} for some constant $C$.
 \bbox

\subsection{The number of excedances in random derangements.} The last example resembles the first two in many aspects, so we will derive the martingale differences and leave it there. The number of excedances in derangements, see the definition in the beginning of Section \ref{sect2}, satisfies the following recursive relation \cite{Z96}.
\begin{equation}
D^{\textnormal{exc}}_{n,k}=kD^{\textnormal{exc}}_{n-1,k}+(n-k)D^{\textnormal{exc}}_{n-1,k-1} + (n-1) D^{\textnormal{exc}}_{n-2,k-1},
\end{equation}
which gives the random process 
\begin{equation*} 
D^{\textnormal{exc}}_{n+2}=\begin{aligned}
   &   
   \begin{cases}
     D^{\textnormal{exc}}_{n}+1,&  \text{with prob. } \frac{(n+1)d_n}{d_{n+2}},\\[5pt]
       D^{\textnormal{exc}}_{n+1},&  \text{with prob. } \frac{(n+1)d_{n+1}}{d_{n+2}}\frac{D^{\textnormal{exc}}_{n+1}}{(n+1)},\\[5pt]
      D^{\textnormal{exc}}_{n+1}+1,&  \text{with prob. } \frac{(n+1)d_{n+1}}{d_{n+2}}\frac{n+1-D^{\textnormal{exc}}_{n+1}}{(n+1)}.
    \end{cases}
\end{aligned}
\end{equation*}
Let $Z_n:=(n-1)(D_{n}^{exc}-\mu_n).$  The differences for one-jumps are
\begin{equation*}
Z_i-Z_{i-1}=X_{i,1}+\alpha_{i,1} 
\end{equation*}
where 
\begin{equation*}
X_{i,1}=\begin{aligned}
   &   
   \begin{cases}
     D^{\textnormal{exc}}_{i-1}-(i-1),&  \text{with prob. } \frac{D^{\textnormal{exc}}_{i-1}}{(i-1)},\\[5pt]
     D^{\textnormal{exc}}_{i-1},&  \text{with prob. } \frac{i-1-D^{\textnormal{exc}}_{i-1}}{(i-1)}.
    \end{cases}
\end{aligned}
\end{equation*}
and $\alpha_{i,1}=i-1-(i-1)\mu_{i}+(i-2)\mu_{i-1}.$ For two-jumps, we have
\begin{equation*} 
Z_i-Z_{i-2}=X_{i,2}+ \alpha_{i,2}
\end{equation*}
where $X_{i,2}=2(D^{\textnormal{exc}}_{i-2}-\mu_{i-2})$ and $\alpha_{i,2}=2i-2-(i-1)\mu_{i}+(i-1)\mu_{i-2}.$ So we can write $Z_n$ conditioned on some $\mathbf{a} \in \mathcal{C}_{\leq 2}(n)$ as the sum of a deterministic term in addition to a sum of martingale differences,
\begin{equation*}
\mathbf{E}[Z_n | C_{\leq 2}(n)=\mathbf{a}]=\sum_i \left( X_{\mathbf{a_i}}+ \alpha_{\mathbf{a_i}} \right).
\end{equation*}
A proof of the asymptotic normality of $D_n^{\textnormal{exc}}$ is given in \cite{C02} along with its first two moments. The rate of convergence in the limit can be studied by Theorem \ref{convthm} combined with Rosenthal's inequality and the bounds derived for the random variables defined in Section \ref{sect:app}. 

\section{Higher order recurrences and other statistics} \label{sect5}

Although all the examples that have been covered are derived from triangular arrays obeying a second-order recurrence relation, the same methods can be applied to statistics of higher order recurrence relations. For a recurrence relation of order $s$, we can extend the update rule \eqref{rule} as follows. Given that $\zeta^{(i)}$ is the $i$th coordinate of $\zeta,$
\begin{equation*}
\zeta_{n+1}^{(i)}=\zeta_{n}^{(i+1)} \textnormal{  for all } 2\leq i\leq s, \textnormal{ and }
\end{equation*}
\begin{equation*}
\zeta_{n+1}^{(1)}= \left\lbrace  \zeta_n^{(j)}+ X_{n,s+1-j}\text{ with probability } q_j, \, 1 \leq j \leq s \right\rbrace \quad \textnormal{where} \quad \sum_{j=1}^s q_s =1.
\end{equation*}
So the decomposition can be achieved over compositions of summands of size at most the order of the recurrence. The generating functions for Eulerian statistics that would satisfy higher order recurrence relations, which counts the descents in a given conjugacy class, are obtained in \cite{Fu}. The application of Zeilberger's algorithm (see Section 4.3 of \cite{Ai} for a comparative account of it) in \cite{GZ} suggests that the order of the recurrence is related to the length of the maximum cycle in the conjugacy class or in the union of the conjugacy classes.  

Another type of Eulerian statistic that can possibly be studied by the methods in this paper, is Eulerian-Fibonacci numbers, defined by Carlitz in \cite{C78}. They satisfy the recurrence relation
\begin{equation*}
F_{n,k}= k F_{n-1,k} + (n-k+1)F_{n-1,k-1} + F_{n-2,k} - 2 F_{n-2,k-1}+ F_{n-2,k-2}.
\end{equation*}
The contribution of the second term is independent of $k,$ but the difficulty is due to the middle term of the second-degree update, which is negative. The asymptotic normality of these numbers is argued to follow from the method of moments in \cite{H20}, along with the study of many other examples of Eulerian statistics described in various recursive forms. 
 
Another direction to follow is the study of Eulerian statistics in different Coxeter groups. As observed in \cite{Oz}, the martingale methods can be extended to Coxeter groups of type B thanks to the equivalent recurrence relations of the symmetric group up to a constant factor. However, the Coxeter group of type D does not allow a similar extension of the methods. Even the descents in random elements of the group with no restriction have complicated recurrence relations, which can be found in \cite{C08}, that do not immediately yield jump probabilities. Yet, the asymptotic normality is already known as their generating function is real-rooted \cite{SV}. 

Some other examples that are not Eulerian but still satisfy a simple recursive relation include variations of Fibonacci permutations. See \cite{K19} for numerous examples of those numbers. The techniques employed in Section \ref{sect4} do not seem to apply to most of them, but the rudimentary methods in Section $3$ can be improved to approximate their moments at least. Another reason to study them is their potential to reveal useful identities. An illustrative example of a Fibonacci-variant statistic is introduced in \cite{H76} by the following recurrence relation
\begin{equation*}
    H_{n,k} = H_{n-1,k} + H_{n-2,k}= H_{n-1,k-1} + H_{n-2, k-2}
\end{equation*}
and the initial condition $H_{1,0}=1.$ It defines a triangular array, which in fact has convolutions of Fibonacci numbers as its rows, i.e., $H_{n,k}= f_{k+1}f_{n-k+1}.$ Diaconis in \cite{D18} defined a Markov chain to study the distribution of bits of binary strings with no two consecutive ones, and the probability of the $k$th bit to be zero is indeed proportional to $H_{n,k}.$ The result in the paper implies 
\begin{equation*}
\mathbf{P}(H_n=k)=C \left(1+ \mathcal{O}\left(\varphi^{-2k}\right)+ \mathcal{O}\left(\varphi^{-2(n-k)}\right) \right).
\end{equation*}
This suggests that it converges to the uniform distribution, which is not decomposable, other than the spikes on both its ends.

The final remark is about the technical aspects of the random process defined in Section \ref{sect3}, which can be viewed as an autoregressive model in the form
\begin{equation*}
Y_n = A_n Y_{n-1}+ B_n
\end{equation*}
where $Y_i$ and $B_i$ are $d$-dimensional vectors and $A_i$ is a $d\times d$ matrix. In particular, \eqref{rule} can be written as
\begin{equation*}
\zeta_{n}^T=\begin{bmatrix}
1-q_n & q_n \\ 1 & 0 
\end{bmatrix} \zeta_{n-1}^T + \begin{bmatrix}
X_{n,2} \\ X_{n,1}
\end{bmatrix}.
\end{equation*}
The conditions for the existence of a stationary solution analogous to the formula \eqref{inhom} are given in \cite{BP92}. In \cite{K73}, a central limit theorem is shown provided that $A_n$ and $B_n$ are independent and identically distributed. The case with martingale differences is addressed in \cite{AK92} with a multivariate limit theorem, whose assumptions include that $A_n$ is constant and the covariance matrix for $B_n$ converges in probability. Although none of the two holds in our case, it could be possible to weaken them and make inferences about the limiting behavior of processes derived from Eulerian statistics without recourse to decompositions. 

\bibliographystyle{alpha}
\bibliography{involutions}

\end{document}